\RequirePackage{fix-cm}
\documentclass{article}                     
\usepackage{graphicx}
\usepackage[margin=1.0in]{geometry}

\usepackage{amsmath}
\usepackage{algorithm}
\usepackage{epsfig}
\usepackage{subfigure}
\usepackage{graphics}
\usepackage{color,ulem}
\usepackage{amssymb}
\usepackage{algorithm,hyperref}
\usepackage[noend]{algpseudocode}
\usepackage[numbers,sort]{natbib}
\usepackage{authblk}

\usepackage{graphicx,booktabs,array}
\newcolumntype{C}{>{\centering\arraybackslash}m{10em}}

\begin{document}

\title{Hierarchical model reduction driven by machine learning \\for parametric advection-diffusion-reaction problems \\in the presence of noisy data
}

\author[*]{Massimiliano Lupo Pasini }
\author[**]{Simona Perotto }
\affil[*]{Computational Sciences and Engineering Division\\ Oak Ridge National Laboratory\\ 
    1 Bethel Valley Road, Oak Ridge, TN, USA, 37831}
\affil[**]{MOX -- 
	Dipartimento di Matematica\\
	Politecnico di Milano\\
	Piazza L. da Vinci, 32, I-20133 Milano, Italy\\}

\renewcommand\Authands{ and }


\date{}

\newtheorem{remark}{Remark}

\maketitle

\footnotetext[1]{\noindent $^*$ lupopasinim@ornl.gov \\ $^{**}$ simona.perotto@polimi.it}

\begin{abstract}
We propose a new approach to 
generate a reliable reduced model for a parametric elliptic problem, in the presence of noisy data. The reference model reduction procedure is the directional HiPOD method, which 
combines Hierarchical Model reduction with a standard Proper Orthogonal Decomposition, according to an offline/online paradigm.
In this paper we show that directional HiPOD looses in terms of accuracy when problem data are affected by noise. This is due to the interpolation driving the online phase, since it replicates, by definition, the noise trend. To overcome this limit, we replace interpolation with Machine Learning fitting models which better discriminate relevant physical features in the data from irrelevant unstructured noise. 
The numerical assessment, although preliminary, confirms the potentialities of the new approach.
\end{abstract}

{\footnotesize \noindent This manuscript has been authored in part by UT-Battelle, LLC, under contract DE-AC05-00OR22725 with the US Department of Energy (DOE). The US government retains and the publisher, by accepting the article for publication, acknowledges that the US government retains a nonexclusive, paid-up, irrevocable, worldwide license to publish or reproduce the published form of this manuscript, or allow others to do so, for US government purposes. DOE will provide public access to these results of federally sponsored research in accordance with the DOE Public Access Plan (\url{http://energy.gov/downloads/doe-public-access-plan}).
}

\section{Introduction and motivations}\label{sec1}
Numerical methods for solving parametric partial differential equations are relevant for all the engineering applications
that can be framed 
into a multi-query or a real-time context. Recurrent instances 
range from the estimation of the parameters which govern physical phenomena (e.g., in biomedical engineering, where the known velocity of the blood and the Navier-Stokes equations for incompressible fluids are combined to draw conclusions about the blood viscosity inside a cardiovascular vessel) to the solution of optimal control problems (e.g., in environmental engineering, where experimental measurements of the concentration of a pollutant in the water and the Navier-Stokes equations for incompressible fluids with a parametrized forcing term are combined to compute the maximum amount of pollutant that can be released without compromising the ecosystem inside a river).

The curse of dimensionality characterizing full order models to simulate this type of physical systems raised the necessity to propose specific numerical methods in order to sustain the computational cost.
Reduced order models have garnered interest in the scientific computing community as an effective way to compress complex partial differential models in a projection subspace, and thus make it more computationally convenient to solve \cite{Audouze2009,Benner11,Benner15,ghattas_willcox_2021,HesthavenRozzaStamm16,QuarteroniManzoniNegri16}.

Here, we focus on a class of reduced order models conceived to describe flows in pipes or, more in general, phenomena with a privileged dynamics aligned with the centerline of the pipe, which may be locally modified by secondary dynamics evolving along the transverse sections. We are referring to the Hierarchical Model (HiMod)
reduction~\cite{ErnPerottoVeneziani08,PerottoErnVeneziani10,Perotto14b,PerottoZilio13}, and, in particular, to a parametric counterpart of such an approach, known as HiPOD~\cite{Barolietal17,lupopasini2022}. HiPOD offers a possible remedy to the well-known bottleneck of an offline/online paradigm, i.e., the computational burden characterizing the offline phase. The idea is to replace the “truth” model in the offline phase of the Proper Orthogonal Decomposition (POD)~\cite{KahlbacherVolkwein07,Kerschen2005,KunishVolkwein02,Volkweinn13} with
a HiMod discretization, namely, more in general, with a
reduced-order model characterized by a high accuracy and a contained computational demand. Then, the online phase recovers the HiMod approximation for a not yet sampled value of the selected parameter, after solving a problem of a very small dimensionality. The computational advantages led by a HiPOD approximation have been numerically investigated both on scalar and vector problems~\cite{lupopasini2022,Barolietal17}.\\
In~\cite{lupopasini2022} two HiPOD model reduction procedures are presented. Here, we consider the directional approach, that combines HiMod and POD by exploiting the decoupling between leading and secondary dynamics at the basis of a HiMod formulation.

In this work, we show that the presence of noise in the data can make directional HiPOD unstable, due to an interpolation over the parametric
space to estimate the online solution. Actually, interpolation reproduces, instead of filtering out, the noise affecting data. 
To overcome this issue, we propose to replace interpolation with 
Machine Learning (ML) fitting models which better discriminate relevant physical features in the data from irrelevant unstructured noise, thus allowing the model to retain and
leverage the former and discard the latter. The numerical assessment we carried out shows that the ML methodology reduces the $L^2$- and the $H^1$-norm of the HiPOD relative error up to an order of magnitude with respect to state-of-the-art interpolation techniques, thus confirming the strong potentiality of the proposed approach.


The paper is structured as follows. Section~\ref{sec_BG} addresses the mathematical background exploited to upgrade the directional HiPOD procedure in~\cite{lupopasini2022}, namely the basics to perform a HiPOD model reduction and some machine learning regression models.
Section~\ref{secnewmethod} introduces the new HiPOD approach and analyzes the effect of the noise propagation onto the response matrix. In Section~\ref{num_sec}, we carry out a numerical assessment of the proposed methodology, by comparing the state-of-the-art HiPOD approach with the new ML version, both in terms of accuracy and robustness to the level of noise.
Finally, some conclusions and perspective are supplied in the last section.

\section{Mathematical background}\label{sec_BG}
HiPOD model reduction and Machine Learning (ML) models for a regression analysis represent the main methodological tools supporting the new approach proposed in this paper.
The next sections are devoted to introduce such tools, for the reader completeness.\\
In particular, we choose as a reference problem a parametrized elliptic Partial Differential Equation
(PDE), defined on a domain $\Omega\subset \mathbb R^2$, whose weak form is 
\begin{equation}\label{parametrized_referenc}
\mbox{find\ } u(\alpha)\in V\quad \mbox{s.t.} \quad 
a(u(\alpha), v; \alpha) = f(v; \alpha) \quad \forall v\in V,
\end{equation}
where $\alpha\in {\mathcal P}\subset \mathbb{R}^p$ is the selected parameter varying in the set ${\mathcal P}$ of the admissible values;
$a(\cdot, \cdot; \alpha): V \times V \times {\mathcal P}\rightarrow \mathbb{R}$ 
and $f(\cdot; \alpha): V \times {\mathcal P} \rightarrow \mathbb{R}$ denote the parametrized bilinear 
and linear forms characterizing the differential problem at hand, the linearity being
meant with respect to the variables different from $\alpha$, and with 
$V\subseteq H^1(\Omega)$ a suitable Hilbert space depending on the specific PDE problem as well as on the assigned boundary conditions, 
standard notation being adopted for function spaces~\cite{ErnGuermond}.
To simplify the exposition, we focus here on the case of a standard scalar linear advection-diffusion-reaction (ADR)
problem, completed with full homogeneous Dirichlet boundary conditions (see Section~\ref{num_sec} for a more general setting), so that
the bilinear and the linear forms in \eqref{parametrized_referenc} are 
\begin{equation}\label{forms}
a(w, z; \alpha) = \displaystyle \int_\Omega \mu \nabla w \cdot  \nabla z \, d\Omega
+  \int_\Omega \big( {\bf b} \cdot \nabla w + \sigma w \big) z  \, d\Omega, \quad
f(z; \alpha) = \displaystyle \int_\Omega fz  \, d\Omega,
\end{equation}
with $w$, $z\in V=H^1_0(\Omega)$, and where parameter $\alpha$ coincides with some of the problem data, i.e.,
the viscosity $\mu$, the advective field ${\bf b}=[b_1, b_2]^T$, the reaction $\sigma$, the source term $f$ (or a boundary value when non homogeneous or more general boundary conditions are assigned).\\
Suitable assumptions are advanced  on the problem data in order to guarantee the well-posedness of formulation \eqref{parametrized_referenc}, for any $\alpha \in {\mathcal P}$.
Finally, we conjecture an affine parameter dependence~\cite{HesthavenRozzaStamm16,QuarteroniManzoniNegri16}.

\subsection{The HiPOD approach}
The HiPOD method provides the parametric counterpart of a Hierarchically Model (HiMod) reduction, by properly combining a
HiMod discretization with the standard Proper Orthogonal Decomposition (POD)~\cite{KahlbacherVolkwein07,Kerschen2005,KunishVolkwein02,Volkweinn13}.
In the next section we recap the main features characterizing a HiMod reduction, being instrumental in the setting of HiPOD procedures.

\subsubsection{Hierachical model reduction}
HiMod reduction proved to be an ideal method for the modeling of scenarios where a dominant direction is evident in the global dynamics of the considered phenomenon. This modeling property is, in general, mirrored by the geometric design of the computational domain, which is assumed to coincide with a pipe where the centerline is parallel to the dominant dynamics. 
Thus, $\Omega$, is identified by
the fiber bundle
$\bigcup_{x \in \Omega_{1D}} \{x\} \times \Sigma_x $, where $\Omega_{1D}$ is the one-dimensional (1D) supporting fiber paraller to the leading dynamics, while $\Sigma_x$ denotes the transverse section at 
the generic point $x$ along $\Omega_{1D}$~\cite{ErnPerottoVeneziani08,PerottoErnVeneziani10,PerottoZilio13,Perotto14b}.
For simplicity, we focus on rectilinear domains, so that 
$\Omega_{1D}\equiv(a,b)\subset
\mathbb{R}$ (we refer the interested reader to~\cite{Perotto14a,PerottoRealiRusconiVeneziani17,PerottoBarbosa20} for pipes with a bent centerline).

HiMod reduction exploits the geometric requirement on $\Omega$ to discretize a PDE problem, so that the leading and the transverse dynamics are approximated with different methods, according to a separation of variable criterion. In this work, we adopt the approach proposed in~\cite{AlettiPerottoVeneziani18} where finite elements model the main dynamics, while a basis of customized modal functions describes the dynamics parallel to the  transverse sections of the pipe. This decoupling of leading and secondary dynamics allows us to commute the full problem into a system of coupled 1D problems, independently of the original model dimensionality, with a consequent considerable benefit in terms of computational effort~\cite{AlettiPerottoVeneziani18,Alvarezetal17,GuzzettiPerottoVeneziani18,PerottoBarbosa20}.\\
Additionally, as it is recurrent in several well-known contexts, computations are performed in a reference domain $\widehat{\Omega}$ (where, e.g., constants can be explicitly computed) and successively moved to the physical domain $\Omega$. This is performed by means of a map $\Psi: \Omega\rightarrow\widehat{\Omega}$ which is assumed to be differentiable with respect to both the independent variables $x$ and $y$. The domain $\widehat{\Omega}$ shares a fiber structure as $\Omega$, being $\widehat \Omega=\Omega_{1D} \times \widehat{\Sigma}$, with $\widehat{\Sigma}$ the reference transverse fiber and where the supporting fiber is the same as for $\Omega$. In particular, for any point
${\bf z}=(x, {\bf y})\in \Omega$,
there exists a point $\widehat {\bf z}=(\widehat x, \widehat {\bf y})\in \widehat{\Omega}$, such that
$\widehat {\bf z}=\Psi({\bf z})$,
with $\widehat x\equiv x$ and $\widehat {\bf y}=\psi_x({\bf y})$,  
$\psi_x: \Sigma_x \rightarrow \widehat{\Sigma}$ being the map between the generic and the reference transverse fiber.
Hereafter, we assume $\psi_x$ to be a $C^1$-diffeomorphism, for all $x\in \Omega_{1D}$ (more details about maps $\Psi$ and $\psi_x$ are available in~\cite{PerottoErnVeneziani10}).

The separation of variables combined with the mapping to $\widehat \Omega$ leads us to define the parametric HiMod reduced space 
\begin{equation}\label{spazioR}
V_{m}(\alpha)=\Big\{  v_{m}(x, {\bf y}; \alpha) = \displaystyle \sum_{k=1}^{m}
\sum_{j=1}^{N_h}  {\tilde v}_{k, j}^\alpha \vartheta_j(x) \varphi_k(\psi_x({\bf y})),\ (x, {\bf y}; \alpha)\in \Omega_{1D} \times\Sigma_x\times {\mathcal P}
 \Big\},
\end{equation}
where $m\in  \mathbb{N}^+$ is the modal index setting the level of detail of the HiMod approximation
in the hierarchy; ${\mathcal B_1}=\{\vartheta_j\}_{j=1}^{N_h}$ is a basis for the 1D space, $V_{1D}\subset H^1_0(\Omega_{1D})$, of the finite element functions associated with the supporting fiber $\Omega_{1D}$ and vanishing at $a$ and $b$, with 
${\rm dim}(V_{1D})=N_h<+\infty$; 
${\mathcal B_2}=\{ \varphi_k\}_{k\in \mathbb{N}^+}$ denotes the basis of modal 
functions defined on the reference transverse fiber $\widehat \Sigma$,  orthonormal with respect to the $L^2(\widehat \Sigma)$-scalar product
and including, in an essential way~\cite{AlettiPerottoVeneziani18}, the data assigned on the lateral boundary $\Gamma_{L}=\cup_{x\in \Omega_{1D}} \partial \Sigma_x$ of $\Omega$ (the reader interested to different possible choices both for basis ${\mathcal B_1}$ and ${\mathcal B_2}$ may refer, e.g., to~\cite{ErnPerottoVeneziani08,PerottoErnVeneziani10,GuzzettiPerottoVeneziani18,PerottoRealiRusconiVeneziani17,PerottoBarbosa20}).
In particular, function ${\tilde v}_{k}(x; \alpha) 
= \sum_{j=1}^{N_h} {\tilde v}_{k, j}^\alpha \vartheta_j(x) \in V_{1D}$ identifies the frequency coefficient
associated with the $k$-th modal function $\varphi_k$.

As far as the modal index $m$ is concerned, it may be assigned thanks to a trial-and-error procedure (see, e.g.,~\cite{ErnPerottoVeneziani08,PerottoErnVeneziani10}) or starting from some preliminary (geometric or physic) information about the problem at hand (as, e.g., in~\cite{GuzzettiPerottoVeneziani18}) or via an automatic selection based on
an a posteriori modeling error analysis (we refer, e.g., to~\cite{PerottoVeneziani14,PerottoZilio15}). We adopt the same value for $m$ in the whole $\Omega$, although the modal index may be locally varied along the domain to match possible heterogeneities of the solution (see~\cite{PerottoErnVeneziani10,PerottoZilio13,PerottoVeneziani14,Perotto14b,PerottoZilio15} for further details).

Thus, the HiMod approximation to problem \eqref{parametrized_referenc}, associated with the modal coefficient $m$, reads
\begin{equation}\label{parametrized_HiMod_reference}
\mbox{find\ } u_{m}(\alpha)\in V_{m}(\alpha)\quad \mbox{s.t.} \quad 
a(u_{m}(\alpha), v_{m}; \alpha) = f(v_{m}; \alpha) \quad \forall v_{m}\in V_{m}(\alpha),
\end{equation}
with $u_{m}(\alpha)=u_{m}(x, {\bf y}; \alpha)$.
To ensure the well-posedness of formulation \eqref{parametrized_HiMod_reference}, we endow the space $V_m(\alpha)$ both with a conformity and a spectral approximability
assumption, while the convergence
of the HiMod approximation $u_{m}(\alpha)$ to the full solution $u(\alpha)$ in \eqref{parametrized_referenc} is guaranteed by introducing a standard density assumption on the discrete space $V_{1D}$.

After applying the HiMod representation in \eqref{spazioR} to $u_{m}(\alpha)$ in \eqref{parametrized_HiMod_reference}, and choosing the test function $v_{m}$ as $\vartheta_t \varphi_q$, for $t=1, \ldots, N_h$ and $q=1, \ldots, m$, the HiMod formulation turns into the system of $m$ coupled 1D problems, 
\begin{equation}\label{HiModsystem}
A_m(\alpha){\bf u}_m(\alpha)={\bf f}_m(\alpha),
\end{equation}
with $A_m(\alpha)\in \mathbb{R}^{mN_h\times mN_h}$ denoting the HiMod stiffness matrix, ${\bf f}_m(\alpha)\in \mathbb{R}^{mN_h}$ representing the HiMod right-hand side, and where vector
\begin{equation}\label{HiModAPP}
{\bf u}_m(\alpha)=\big[ {\tilde u}_{1, 1}^\alpha, \ldots, {\tilde u}_{1, N_h}^\alpha, {\tilde u}_{2, 1}^\alpha, \ldots, {\tilde u}_{2, N_h}^\alpha, \ldots, {\tilde u}_{m, 1}^\alpha, \ldots, {\tilde u}_{m, N_h}^\alpha\big]^T\in \mathbb{R}^{mN_h}
\end{equation}
collects the modal coefficients $\{ {\tilde u}_{k, j}^\alpha \}_{k=1, j=1}^{m, N_h}$, namely the actual unknowns of the HiMod discretization
\begin{equation}\label{HiModexp}
u_{m}(x, {\bf y}; \alpha) = \displaystyle \sum_{k=1}^{m} 
 \sum_{j=1}^{N_h}  {\tilde u}_{k, j}^\alpha \vartheta_j(x) \varphi_k(\psi_x({\bf y}))
\end{equation}
\cite{ErnPerottoVeneziani08,PerottoErnVeneziani10}. It has been numerically checked that, when the mainstream dominates the transverse dynamics (i.e.,
for small values of $m$), the HiMod approximation demands a considerably lower computational effort if compared with a standard (e.g., finite element) discretization of problem 
\eqref{parametrized_referenc}, without giving up the accuracy of the simulation~\cite{AlettiPerottoVeneziani18,Alvarezetal17,GuzzettiPerottoVeneziani18,PerottoBarbosa20}.

\subsubsection{Directional HiPOD reduction}\label{DHiPOD_sec}
A HiPOD approach offers a new way to yield a HiMod approximation, as an alternative to the resolution of the HiMod system \eqref{HiModsystem}.
According to a data-driven procedure,
the idea is to replace the HiMod discretization in \eqref{HiModexp}
with a surrogate solution obtained
by resorting to a POD reduced basis generated by HiMod approximations.
A standard offline/online paradigm drives the computation of the HiPOD approximation. In particular, 
during the offline phase, we
compute the HiMod solution in \eqref{HiModsystem} for different choices of $\alpha$, to extract the POD basis.
In the online step, such a basis is employed to approximate the HiMod solution in \eqref{HiModsystem} associated with a value of the parameter not sampled during the offline phase.\\
As shown in~\cite{Barolietal17,lupopasini2022}, a HiPOD procedure considerably lowers computational costs, without compromising the quality of the reduced solution.

Two HiPOD approaches have been explored so far. The basic method 
coincides with a straightforward application of a projection-based POD to HiMod solutions~\cite{Barolietal17}. An advanced procedure, referred to as directional HiPOD, takes advantage of the separation of variables implied by a HiMod discretization~\cite{lupopasini2022}, 
the SVD being used to remove the redundancy along the main stream and the transverse direction, separately.
In addition, the online phase is carried out by interpolation instead of projection, thus relieving us from assembling
the HiMod stiffness matrix and the right-hand side associated with the online parameter, as expected by the basic HiPOD approach. 

Below, we detail the directional HiPOD method, since instrumental to the new approach proposed in Section~\ref{secnewmethod}. \\
The goal of the offline phase is to identify the POD basis to be used in an online mode. To this aim, we build the response matrix by collecting the HiMod solution to
problem \eqref{parametrized_referenc} for $p$ different values, $\alpha_i$,
of the parameter $\alpha$, with $i=1, \ldots, p$. In more detail, the generic HiMod solution $u_{m}(x, {\bf y}; \alpha_i)$ is identified with the corresponding modal coefficients, $\{ {\tilde u}_{k, j}^{\alpha_i} \}_{k=1, j=1}^{m, N_h}$, collected by mode into the $m$ vectors
\begin{equation}\label{vectU}
{\bf U}^k(\alpha_i)=[\tilde u_{k, 1}^{\alpha_i}, \tilde u_{k, 2}^{\alpha_i},
\ldots, \tilde u_{k, N_h}^{\alpha_i}]^T\in \mathbb{R}^{N_h} \quad k=1, \ldots, m,
\end{equation}
so that the response matrix is assembled as 
\begin{equation}
\begin{array}{rcl}
U&=&\big[ {\bf U}^1(\alpha_1) \cdots {\bf U}^m(\alpha_1)  \lvert  {\bf U}^1(\alpha_2) \cdots {\bf U}^m(\alpha_2)  \lvert  \cdots \cdots \cdots
\lvert \ {\bf U}^1(\alpha_p) \cdots {\bf U}^m(\alpha_p)\big]\qquad \qquad \\[4mm]
&=& \left[
\begin{array}{ccc|ccc|cc|ccc}
\tilde u_{1, 1}^{\alpha_1} & \cdots & \tilde u_{m, 1}^{\alpha_1} & \tilde u_{1, 1}^{\alpha_2} & \cdots & \tilde u_{m, 1}^{\alpha_2}  & 
\cdots & \cdots &  \tilde u_{1, 1}^{\alpha_p} & \cdots & \tilde u_{m, 1}^{\alpha_p} \\
\tilde u_{1, 2}^{\alpha_1} & \cdots & \tilde u_{m, 2}^{\alpha_1} & \tilde u_{1, 2}^{\alpha_2} & \cdots & \tilde u_{m, 2}^{\alpha_2}  & 
\cdots & \cdots & \tilde u_{1, 2}^{\alpha_p} & \cdots & \tilde u_{m, 2}^{\alpha_p}\\
\vdots & \vdots & \vdots & \vdots & \vdots & \vdots & \vdots & \vdots & \vdots & \vdots & \vdots \\
\tilde u_{1, N_h}^{\alpha_1} & \cdots & \tilde u_{m, N_h}^{\alpha_1} & \tilde u_{1, N_h}^{\alpha_2} & \cdots & \tilde u_{m, N_h}^{\alpha_2}  & 
\cdots & \cdots & \tilde u_{1, N_h}^{\alpha_p} & \cdots & \tilde u_{m, N_h}^{\alpha_p}
\end{array}
\right].
\end{array}
\label{response_matrix}
\end{equation}
The blocks are associated with the different parameters, while, for each block, columns run over modes, 
rows over the finite element nodes.
In order to extract the POD basis, 
we apply the Singular Value Decomposition (SVD)~\cite{GolubVanLoan13} to matrix $U$, so that
\begin{equation}\label{Ufac}
U=\Xi \Lambda  K^T, 
\end{equation}
with $\Xi \in \mathbb{R}^{N_h \times N_h}$ and $K \in \mathbb{R}^{(mp)\times (mp)}$ unitary matrices, and
$\Lambda \in \mathbb{R}^{N_h \times (mp)}$ a pseudo-diagonal matrix. 
Each column in $U$ can be expanded in terms of the left singular vectors $\{ {\boldsymbol \xi}_j\}_{j=1}^{N_h}$ of $U$ since constituting an orthogonal basis for $\mathbb{R}^{N_h}$, so that
\begin{equation}\label{vectorU}
{\bf U}^k(\alpha_i) = \displaystyle \sum_{j=1}^{N_h} T_j^k(\alpha_i) {\boldsymbol \xi}_j \quad k=1, \ldots, m, \ i=1, \ldots, p.
\end{equation}
It is customary to select the first, say $L$ with $L\le N_h$, most meaningful singular vectors of $U$ to identify the reduced POD space $V_{\rm POD, 1}^L= \text{span} \{{\boldsymbol \xi}_1, \dots, {\boldsymbol \xi}_{L}\}$, with $\dim(V_{\rm POD, 1}^L)=L$. As a consequence, the vectors in \eqref{vectorU} are approximated as
\begin{equation}\label{vectorUap}
{\bf U}^k(\alpha_i) \cong  \displaystyle \sum_{j=1}^{L} T_j^k(\alpha_i) {\boldsymbol \xi}_j \quad k=1, \ldots, m, \ i=1, \ldots, p,
\end{equation}
the equality being ensured for $L=N_h$. The directional HiPOD reduction involves a reorganization of coefficients $\{T_j^k(\alpha_i)\}$ first into the vectors ${\bf T}_j(\alpha_i)=[T_j^1(\alpha_i), \ldots, T_j^m(\alpha_i)]^T\in \mathbb{R}^m$ with $i=1, \ldots, p$, and successively into the matrices
$$
S_j=[{\bf T}_j(\alpha_1), \ldots, {\bf T}_j(\alpha_p)]=\left[
\begin{array}{ccc}
T_j^1(\alpha_1) & \ldots & T_j^1(\alpha_p)\\
\vdots & & \vdots\\
T_j^m(\alpha_1) & \ldots & T_j^m(\alpha_p)
\end{array}
\right]\in \mathbb{R}^{m\times p},
$$
with $j=1, \ldots, L$, in order to associate a reduced POD basis with each index $j$. To this aim, we factorize the matrices $S_j$ via SVD, so that
\begin{equation}\label{svd_Sj}
S_j=R_j D_j P_j^T,
\end{equation}
with $R_j \in \mathbb{R}^{m\times m}$ and $P_j\in \mathbb{R}^{p\times p}$ unitary matrices, and $D_j\in \mathbb{R}^{m\times p}$ 
the pseudo-diagonal matrix collecting the singular values of $S_j$.
It follows that each column ${\bf T}_j(\alpha_i)$ of $S_j$
can be approximated by resorting to the POD orthogonal basis $\{ {\bf r}_j^k \}_{k=1}^{\mu_j}$, with $\mu_j\le m$, constituted by the most significant $\mu_j$
left singular vectors of $S_j$, namely
\begin{equation}\label{Tcoefficients}
{\bf T}_j(\alpha_i) \cong \displaystyle \sum_{k=1}^{\mu_j} Q_j^k(\alpha_i) {\bf r}_j^k \quad j=1, \ldots, L,\ i=1, \ldots, p.
\end{equation}
Thus, a POD space $V_{{\rm POD, 2}, j}^{\mu_j}= \text{span} \{{\bf r}_j^1, \dots, {\bf r}_j^{\mu_j}\}$ can be defined for each $j$, with $\dim(V_{{\rm POD, 2}, j}^{\mu_j})$ $=\mu_j$.\\
The offline phase ends with the overall generation of $(L+1)$ POD reduced bases
to be used in the online phase in order to predict the HiMod approximation to problem \eqref{parametrized_referenc} for a value, $\alpha^*$,  of the parameter, such that $\alpha^*\ne \alpha_i$ for $i=1, \ldots, p$. 
The idea is to go backward through the directional procedure, starting from 
the coefficients $Q_j^k(\alpha^*)$
in \eqref{Tcoefficients}, with $j=1, \ldots, L$, $k=1, \ldots, \mu_j$, which are approximated by interpolating the 
(known) values  $Q_j^k(\alpha_i)$ for $i=1, \ldots, p$. Coefficients $Q_j^k(\alpha^*)$, together with the $L$ POD bases $\{ {\bf r}_j^k \}_{k=1}^{\mu_j}$, allow us to compute the vectors  
\begin{equation}\label{newQ}
{\bf T}_j(\alpha^*) = [T_j^1(\alpha^*), \ldots, T_j^m(\alpha^*)]^T=
\displaystyle \sum_{k=1}^{\mu_j} Q_j^k(\alpha^*) {\bf r}_j^k \quad j=1, \ldots, L
\end{equation}
in $\mathbb{R}^m$ and, consequently,  by exploiting the POD basis,
$\{ {\boldsymbol \xi}_j\}_{j=1}^{L}$,  generated first, to assemble the $m$ vectors 
\begin{equation}\label{mpod}
{\bf U}_{{\rm HiPOD}}^k(\alpha^*) = [u_{{\rm POD}, k, 1}^{\alpha^*}, \ldots, u_{{\rm POD}, k, N_h}^{\alpha^*}]^T=
\displaystyle \sum_{j=1}^{L} T_j^k(\alpha^*) {\boldsymbol \xi}_j \quad k=1, \ldots, m
\end{equation}
in $\mathbb{R}^{N_h}$, which represent the online counterpart of vectors in \eqref{vectorUap}. The HiMod solution $u_m(\alpha^*)$ can thus be approximated by means of the expansion
\begin{equation}\label{recovery}
u_{{\rm HiPOD}}^{L, M_L}(\alpha^*)=\displaystyle \sum_{k=1}^{m} 
\Big[ \sum_{j=1}^{N_h}  {u}_{{\rm POD}, k, j}^{\alpha^*}\, \vartheta_j(x) \Big] \varphi_k(\psi_x({\bf y})),
\end{equation}
with $M_L=\{\mu_j\}_{j=1}^L$.
In particular, coefficients $u_{{\rm POD}, k, j}^{\alpha^*}$ provide an approximation of the actual coefficient $\tilde u_{k, j}^{\alpha^*}$ in \eqref{vectU} for $\alpha_i=\alpha^*$.

The procedure adopted to predict the coefficients $Q_j^k(\alpha^*)$ plays an important role. When the offline solution data is not affected by noise, interpolation techniques are an effective tool. In~\cite{lupopasini2022}, we assess different intepolations, namely, standard linear interpolation, a piecewise cubic Hermite (PCH) interpolant, and interpolating Radial Basis Functions (RBF), to infer that PCH and RBF interpolants slightly outperform the linear approach.\\
Vice versa, when the offline data is affected by noise, interpolation is not recommended because it does not discriminate between relevant features of the problem and noise, thus not being able to retain the former and discard the latter. This limit  motivated us in the proposal of a new tool, in order to guarantee a reliable HiMod approximation also in the presence of noisy data. 

Finally, to select the dimension of the POD spaces $V_{\rm POD, 1}^L$ and $V_{{\rm POD, 2}, j}^{\mu_j}$, we resort to a control on the variance, i.e., after setting the tolerances $\epsilon_1$ and $\epsilon_2$, with $0\le \epsilon_1, \epsilon_2 \le 1$, we keep the first $L$
left singular vectors ${\boldsymbol \xi}_j$ of $U$ and the first $\mu_j$ left singular vectors ${\bf r}_j^k$ of $S_j$, such that
\begin{equation}\label{tolerances}
\displaystyle \frac{\displaystyle\ \ \sum_{j=1}^L \lambda_j^2\ \ }{\displaystyle\sum_{j=1}^{N_h} \lambda_j^2}\ge \epsilon_1,\quad 
\frac{\ \ \displaystyle \sum_{j=1}^{\mu_j} d_{j, k}^2\ \ }{\displaystyle \sum_{j=1}^{m} d_{j, k}^2}\ge \epsilon_2,
\end{equation}
respectively, with $\lambda_j$ the singular value of $U$ associated with ${\boldsymbol \xi}_j$, for $j=1, \ldots, N_h$, and $d_{j, k}$ the singular value of $S_j$ corresponding to the $k$-th singular vector ${\bf r}_j^k$, with $k=1, \ldots, m$.

\subsection{Machine learning models for regression}
In this section we focus on Machine Learning (ML) models to approximate data distributions.\\ 
The final goal is to replace the interpolation step in the online phase of the directional HiPOD procedure with a regression technique, in order to address situations where the data in the offline phase may be noisy. In particular, to estimate the coefficients $Q_j^k(\alpha^*)$ in \eqref{newQ}, we resort to a ML fitting model since the offline solution is related to parameter $\alpha_i$ by a highly nonlinear relation.

Now, in order to understand how the noise in the data may affect the accuracy of the predictions yielded by a ML fitting model, we have to make some preliminary assumption on the noise properties. 
The most common hypothesis leads us to consider an additive independent identically distributed (i.i.d.) Gaussian noise, $\tilde{\eta}\sim \mathcal{N}(0,\eta)$, in the output, with $\eta>0$, namely, we assume to have
\begin{equation}\label{additive}
    Q_j^k(\alpha_i) + \tilde{\eta}\quad k=1, \ldots, m,\ j=1, \ldots, L,\ i=1, \ldots, p.
\end{equation}
In the next sections, we consider two ML regression models that operate under the additive noise assumption, i.e., the polynomial \cite{Gergonne, Smith, Magee, Fan} and the Gaussian process regression \cite{Dudley, Marcus1, Marcus2, Csato}.

\subsubsection{Polynomial regression}\label{interp_sec}
Polynomial regression is a form of regression analysis where the relationship between the independent and the dependent variables is modeled as a polynomial in the independent variable, of a certain degree $n$. With reference to our specific context, a polynomial regression model of degree $n$ estimates the nonlinear dependence of
$Q_j^k(\alpha_i)$ from the parameter $\alpha_i$
according to the formula
\begin{equation}\label{linear_regression}
 Q_j^k(\alpha_i) =  \sum_{\ell = 0}^{n} \beta_{j, \ell}^k\,  \alpha_i^\ell + \tilde \eta, \quad k=1, \ldots, m,\ j=1, \ldots, L,\ i=1, \ldots, p,
\end{equation}
where $\beta_{j, \ell}^k$ are unknown parameters to be computed in order to optimize the matching between predictions and observations of the dependent variable, while $\tilde{\eta}$ denotes a zero-mean Gaussian noise as in \eqref{additive}. Although relation \eqref{linear_regression} is nonlinear in the independent variable $\alpha_i$, this model is categorized as linear since the regression function is linear in terms of the unknown parameters $\beta_{j, \ell}^k$. \\
With a view to the directional HiPOD procedure,
the quantity
\begin{equation}\label{q*poly}
 \sum_{\ell = 0}^{n} \beta_{j, \ell}^k\,  [\alpha^*]^\ell  \quad k=1, \ldots, m,\ j=1, \ldots, L
\end{equation}
will be used as an estimate for the online coefficient $Q_{j}^k(\alpha^*)$. 

\subsubsection{Gaussian process regression}\label{GP_sec}
The goal of a Gaussian process regression is to determine the best set of random variables that describes the relation between input features and outputs. In the specific setting of interest, we denote the random variables we are looking for by $\mathcal{Q}_{j}^k$, while  $\alpha_i$ and $Q_j^k(\alpha_i)$ represent the inputs and the outputs, respectively.\\
A Gaussian process is fully characterized by its mean function and covariance function. Therefore, Gaussian process regression reduces to calculating the best values for the mean and the covariance functions in order to minimize the mismatch between
outputs ($Q_j^k(\alpha_i)$) and the 
predictions produced by the finite set of random variables selected from the Gaussian process, evaluated at the available inputs ($\alpha_i$). 

Bayesian statistics combines prior Gaussian processes, that retain preliminary knowledge and information from the offline data samples (also known as likelihood), to construct an updated posterior Gaussian process. 
According to the Bayes' rule~\cite{Bayes}, a prior Gaussian process is iteratively updated using the information from the data till it converges to a stationary state defined by a posterior Gaussian process.\\ 
Now, we specify such a workflow onto our setting of interest.
We denote a prior Gaussian process by $\mathcal{Q}_{j,\text{prior}}^k$, so that 
\begin{equation}\label{prior}
     \mathcal{Q}_{j,\text{prior}}^k \sim \mathcal{GP} (q^k_{j,\text{prior}}, z^k_{j,\text{prior}} + \eta)
    \quad k=1, \ldots, m,\ j=1, \ldots, L,
\end{equation}
where $q^k_{j,\text{prior}}=q^k_{j,\text{prior}}(\alpha)$ is the mean function, $z^k_{j,\text{prior}}=z^k_{j,\text{prior}}(\alpha, \tilde \alpha)$ is the covariance function, and $\eta$ is the noise function.
The restriction of the prior Gaussian process in \eqref{prior} at the input points $\alpha_i$, for $i=1,\ldots,p$, is
a multivariate $p$-dimensional Gaussian distribution  $\mathcal{N}(\mathbf{m}_{j, \text{prior}}^k, Z_{j, \text{prior}}^k)$,
where the mean vector, $\mathbf{m}_{j, \text{prior}}^k=[m_{j, \text{prior}, r}^k]\in \mathbb R^p$, and the covariance matrix, $Z_{j, \text{prior}}^k=[Z_{j, \text{prior}, rs}^k]\in \mathbb R^{p\times p}$, are identified by relations
\begin{equation}
    m_{j, \text{prior}, r}^k =
    q^k_{j,\text{prior}}(\alpha_r)=\frac{1}{4}\alpha_r ^ 2 \quad r=1,\ldots, p,
    \label{m_vector}
\end{equation}
\begin{equation}
    Z_{j, \text{prior}, rs}^k =
    z^k_{j,\text{prior}}(\alpha_r, \alpha_s)=\exp\Big[-\frac{1}{2}(\alpha_r-\alpha_s)^2\Big] \quad r, s=1,\ldots, p,
\end{equation}
respectively.
The Gaussian process $\mathcal{Q}_{j,\text{like}}^k$ associated with the data, 
is defined as
\begin{equation}\label{likeeq}
\begin{aligned}
    \mathcal{Q}_{j,\text{like}}^k  \sim \mathcal{GP} (q_{j,\text{like}}^k, z_{j,\text{like}}^k) \quad & k=1, \ldots, m,\ j=1,
    \ldots, L,
\end{aligned}
\end{equation}
with $q^k_{j,\text{like}}$ denoting the mean function and $z^k_{j,\text{like}}$ the covariance function. The Gaussian process  $\mathcal{Q}_{j,\text{like}}^k$ is not fully characterizable using the data, meaning that the mean and the covariance functions cannot be uniquely determined. However, the values attained by $q^k_{j,\text{like}}$ and $z^k_{j,\text{like}}$ at the input points $\alpha_i$, for $i=1,\ldots,p$, are known and correspond to a $p$-dimensional Gaussian distribution $\mathcal{N}(\mathbf{m}^k_{j,\text{like}}, Z^k_{j,\text{like}})$, where the mean vector $\mathbf{m}^k_{j,\text{like}}=[m^k_{j,\text{like}, r}]\in \mathbb R^p$ 
and the covariance matrix $Z^k_{j,\text{like}}=[Z^k_{j, \text{like}, rs}]\in \mathbb R^{p\times p}$ 
are defined by
$$
m^k_{j,\text{like}, r}=Q_j^k(\alpha_r)\quad r=1,\ldots, p,
$$
\begin{equation}
    Z^k_{j, \text{like}, rs} = \exp\Big[-\frac{1}{2}\big(Q_j^k(\alpha_r)-Q_j^k(\alpha_s)\big)^2\Big] \quad r, s=1,\ldots, p.
\end{equation}
Now, the posterior can be used to make predictions for unseen values of $\alpha$ (namely, to predict coefficients $Q_j^k(\alpha^*)$ with $\alpha^*$ the online parameter). 
The joint distribution of the Gaussian processes evaluated at the points $\alpha_i$, for $i=1,\ldots,p$, is a $2p$-dimensional Gaussian distribution
\begin{equation}
 \mathcal{N} 
    \bigg( 
    \begin{bmatrix}
    \mathbf{m}_{j,\text{prior}}^k\\[2mm] \mathbf{m}^k_{j,\text{like}}
    \end{bmatrix}, 
    \begin{bmatrix}
    Z^k_{j,\text{prior}} & Z^k_{j,\text{pl}} \\[2mm] Z^{k}_{j,\text{pl}}& Z^k_{j,\text{like}}
    \end{bmatrix}     
    \bigg),
\end{equation}
where the symmetric matrix $Z^k_{j,\text{pl}}\in \mathbb{R}^{p\times p}$ is the correlation matrix between prior and likelihood.
Using Bayes' Theorem, the posterior $p$-dimensional Gaussian distribution is
\begin{equation}\label{post}
 \mathcal{N} \bigg(\mathbf{m}_{j,\text{prior}}^k
    - Z^{k}_{j,\text{pl}}\, [Z^k_{j,\text{prior}}]^{-1}\mathbf{m}^k_{j,\text{like}}, 
    Z^k_{j,\text{prior}} - Z^k_{j,\text{pl}} [ Z^k_{j,\text{like}}]^{-1}Z^{k}_{j,\text{pl}}\bigg).
\end{equation}
The posterior Gaussian process associated with the distribution in \eqref{post} is consequently given by
\begin{equation}\label{posteriori}
    \mathcal{Q}^k_{j,\text{posterior}} \sim \mathcal{G}\mathcal{P} (q^k_{j,\text{posterior}}, z^k_{j,\text{posterior}}).
\end{equation}
The evaluation of the mean function $q^k_{j,\text{posterior}}$
at the new input parameter value $\alpha^*$ is
\begin{equation}\label{q*_GP}
    q^k_{j,\text{posterior}}(\alpha^*) = q^k_{j,\text{prior}}(\alpha^*) - [Z(\boldsymbol{\alpha}, \alpha^*)]^T\,  [Z_{j, \text{prior}}^k]^{-1}\, {\bf m}^k_{j,\text{like}},
\end{equation}
with $q^k_{j,\text{prior}}(\alpha^*)$ defined as in \eqref{m_vector}, 
$Z(\boldsymbol{\alpha}, \alpha^*)\in \mathbb{R}^p$ denoting the covariance vector between the offline data sampled at $\boldsymbol{\alpha}=[\alpha_1, \ldots, \alpha_p]^T\in \mathbb{R}^p$
and the new sampled parameter $\alpha^*$. 
The covariance function in \eqref{posteriori} evaluated at $\alpha^*$ is 
\begin{equation}
    z^k_{j,\text{posterior}}(\alpha^*) = z^k_{j,\text{prior}}(\alpha^*, \alpha^*) - [Z(\boldsymbol{\alpha}, \alpha^*)]^T [Z_{j, \text{like}}^k]^{-1}\, Z(\boldsymbol{\alpha}, \alpha^*),
\end{equation}
with $z^k_{j,\text{prior}}(\alpha^*, \alpha^*)=1$.
The value $q^k_{j,\text{posterior}}(\alpha^*)$ will be used as an estimate for $Q_{j}^k(\alpha^*)$. 

\section{HiPOD reduction for data affected by Gaussian noise}\label{secnewmethod}
To address situations where the offline data is noisy, we propose here a variant of the directional HiPOD reduction presented in Section~\ref{DHiPOD_sec}.

The idea is very straightforward. In the online phase we replace the initial interpolation step used to estimate  coefficients $Q_j^k(\alpha^*)$ with a regression technique.
In particular, formulas \eqref{q*poly} and \eqref{q*_GP} provide us the desired estimate for the coefficients $Q_j^k(\alpha^*)$, when resorting to a polynomial or to a Gaussian process regression, respectively.
Successively, the online phase is performed exactly as in Section~\ref{DHiPOD_sec}, going through
the reconstructions \eqref{newQ}-\eqref{mpod}, to obtain the final expansion in \eqref{recovery}.\\
The choice for the ML regression models in Sections~\ref{interp_sec}-~\ref{GP_sec} is motivated by the highly nonlinear dependence of the offline HiMod solution onto the offline parameters.

The improvement led by the new HiPOD approach is numerically checked in Section~\ref{num_sec}.
In the next section, we list some numerical quantities that can help us monitor the noise propagation throughout the HiPOD procedure.

\subsection{Noise propagation on the response matrix}
In this section we compare the standard directional HiPOD procedure and the new variant proposed in this paper for different noise levels in the offline data. 
In particular, hereafter, we refer to the standard and to the new directional HiPOD approach as to the interpolation-HiPOD and the regression-HiPOD, respectively.

The data noise affects the HiPOD approximation by perturbing the $(L+1)$ SVD's involved in the offline phase. 
In~\cite{lupopasini2022}, we showed that an inaccurate calculation of the first SVD (i.e., of the SVD of the response matrix $U$) compromises the reliability of the directional interpolation-HiPOD. 
For this reason, here we focus on the effect of the data noise onto the decomposition in \eqref{Ufac}.\\
Let us assume that a Gaussian noise affects the linear form in \eqref{parametrized_referenc} (e.g., by perturbing the source term $f$ or the boundary data of the PDE problem at hand). As a consequence, the HiMod system in \eqref{HiModsystem} is replaced by the perturbed problem 
\begin{equation}\label{perturbed_system}
    A_m(\alpha){\tilde{\bf u}}_m(\alpha)={\bf f}_m(\alpha) + \boldsymbol{\eta}_m,
\end{equation}
where $\boldsymbol{\eta}_m \sim \mathcal{N}(\mathbf{0}_{mN_h}, \eta I_{mN_h})$ identifies the noise, with $\mathbf{0}_{mN_h} \in \mathbb{R}^{mN_h}$ the null vector, $I_{mN_h} \in \mathbb{R}^{{mN_h}}\times \mathbb{R}^{{mN_h}}$ the identity matrix, $\eta >0$
the noise level, and ${\tilde{\bf u}}_m(\alpha)$ denotes the associated noisy HiMod discretization.\\
The solution to \eqref{perturbed_system} can be 
regarded as a perturbation of the HiMod solution, ${\bf u}_m(\alpha)$, in \eqref{HiModsystem} by an additive white noise ($A^{-1}_m(\alpha) \boldsymbol{\eta}_m$), being
\begin{equation}
    \tilde {\bf{u}}_m(\alpha) = A^{-1}_m(\alpha){\bf f}_m(\alpha) + A^{-1}_m(\alpha) \boldsymbol{\eta}_m
    = {\bf u}_m(\alpha) + A^{-1}_m(\alpha) \boldsymbol{\eta}_m.
\end{equation}
To simplify the notation, we define the random variable $\tilde{\boldsymbol{\eta}}_m(\alpha) = A^{-1}_m(\alpha) \boldsymbol{\eta}_m$, with 
\begin{equation}
    \tilde{\boldsymbol{\eta}}_m(\alpha) \sim \mathcal{N}(\mathbf{0}, \eta [A_m(\alpha) A^T_m(\alpha)]^{-1} ),
\end{equation}
so that the solution $\tilde {\bf {u}}_m(\alpha)$ to the perturbed HiMod linear system \eqref{perturbed_system} can be recast as 
\begin{equation}
    \tilde {\bf u}_m(\alpha) = {\bf u}_m(\alpha) + \tilde{\boldsymbol{\eta}}_m(\alpha). 
\end{equation}
This decomposition finds a counterpart when assembling the response matrix, $\tilde U$, associated with the offline noisy data. Indeed, thanks to the additive property assumed for the noise ${\boldsymbol{\eta}}_m$,
matrix $\tilde U$, which collects the perturbed HiMod solutions $\tilde {\bf{u}}_m(\alpha_i)$ in \eqref{perturbed_system}, for the $p$ values, $\alpha_1, \ldots, \alpha_p$, of the parameter $\alpha$, can be conceived as a perturbation of matrix $U$
in \eqref{response_matrix}.
Thus, 
after introducing the noise matrix 
\begin{equation}\label{response_error_matrix}
\begin{array}{rcl}
E&=&[ {\bf E}^1(\alpha_1) \cdots {\bf E}^m(\alpha_1)  \lvert  {\bf E}^1(\alpha_2) \cdots {\bf E}^m(\alpha_2)  \lvert  \cdots \cdots \cdots
\lvert \ {\bf E}^1(\alpha_p) \cdots {\bf E}^m(\alpha_p)] \\[4mm]
&=& \left[
\begin{array}{ccc|ccc|cc|ccc}
\tilde \eta_{1, 1}^{\alpha_1} & \cdots & \tilde \eta_{m, 1}^{\alpha_1} & \tilde \eta_{1, 1}^{\alpha_2} & \cdots & \tilde \eta_{m, 1}^{\alpha_2}  & 
\cdots & \cdots &  \tilde \eta_{1, 1}^{\alpha_p} & \cdots & \tilde \eta_{m, 1}^{\alpha_p} \\[2mm]
\tilde \eta_{1, 2}^{\alpha_1} & \cdots & \tilde \eta_{m, 2}^{\alpha_1} & \tilde \eta_{1, 2}^{\alpha_2} & \cdots & \tilde \eta_{m, 2}^{\alpha_2}  & 
\cdots & \cdots & \tilde \eta_{1, 2}^{\alpha_p} & \cdots & \tilde \eta_{m, 2}^{\alpha_p}\\[2mm]
\vdots & \vdots & \vdots & \vdots & \vdots & \vdots & \vdots & \vdots & \vdots & \vdots & \vdots \\[2mm]
\tilde \eta_{1, N_h}^{\alpha_1} & \cdots & \tilde \eta_{m, N_h}^{\alpha_1} & \tilde \eta_{1, N_h}^{\alpha_2} & \cdots & \tilde \eta_{m, N_h}^{\alpha_2}  & 
\cdots & \cdots & \tilde \eta_{1, N_h}^{\alpha_p} & \cdots & \tilde \eta_{m, N_h}^{\alpha_p}
\end{array}
\right],
\end{array}
\end{equation}
where the vectors  
$$
{\bf E}^k(\alpha_i)=[\tilde \eta_{k, 1}^{\alpha_i}, \tilde \eta_{k, 2}^{\alpha_i},
\ldots, \tilde \eta_{k, N_h}^{\alpha_i}]^T\in \mathbb{R}^{N_h} \quad k=1, \ldots, m,\ i=1, \ldots, p,
$$
gather, by mode, the noise affecting the modal coefficients $\{ {\tilde u}_{k, j}^{\alpha_i} \}_{k=1, j=1}^{m, N_h}$ in \eqref{HiModAPP} for $\alpha=\alpha_i$,
the response matrix $\tilde{U} \in \mathbb{R}^{N_h\times (mp)}$ associated with the noisy HiMod solutions coincides with 
\begin{equation}
    \tilde{U} = U + E,
\end{equation}
with $U$ the response matrix in \eqref{response_matrix}. The SVD 
\begin{equation}\label{svd_perturbed}
\tilde U = \tilde \Xi \tilde\Lambda \tilde K
\end{equation}
will consequently replace the factorization of matrix $U$ in \eqref{Ufac}, with  
$\tilde \Xi\in \mathbb R^{N_h\times N_h}$ and $\tilde K\in \mathbb R^{(mp)\times (mp)}$ the unitary matrices of the left and of the right singular vectors of $\tilde U$, and $\tilde \Lambda\in \mathbb R^{N_h\times (mp)}$ the pseudo-diagonal matrix of the singular values.

\begin{remark}
The injection of a noise in the problem data involved in the bilinear form in \eqref{forms}
would still affect the HiMod solution ${\bf u}_m(\alpha)$, in \eqref{HiModsystem}, but not in an additive fashion. This would unavoidably make the perturbation propagation analysis more complex, and is beyond the purpose of this paper.
\end{remark}

Some results are available in the literature which relate both the singular values and the singular vectors of matrices $U$ and $\tilde U$. 

As for the singular values, we remind the Weyl theorem~\cite{Weyl} and the Mirsky theorem~\cite{Mirsky}. The first result controls the discrepancy between the $i$-th singular value $\lambda_i$ 
of $U$ and the corresponding singular value $\tilde \lambda_i$ of $\tilde U$ in terms of the spectral norm, $\lVert E \rVert_S$, of the noise matrix, being   
\begin{equation}\label{weyl}
\lvert \tilde \lambda_i -\lambda_i \rvert \le \lVert E \rVert_S \quad i=1,\ldots,\min\{N_h, mp\}.
\end{equation}
Mirsky theorem provides an upper bound on the sum of the 
quadratic deviations of values $\tilde \lambda_i$'s with respect to $ \lambda_i$'s in terms of the Frobenius norm, $\lVert E \rVert_F$, of the noise matrix, given by   
\begin{equation}\label{mirsky}
 \sum_{i=1}^{\min\{N_h, mp\}}( \tilde \lambda_i -\lambda_i )^2 \le \lVert E \rVert_F. 
\end{equation}
Of course, inequalities \eqref{weyl} and \eqref{mirsky} are completely useless when norms $\lVert E \rVert_S$ and $\lVert E \rVert_F$ become larger and larger.

The reference result on the singular vectors is represented by the generalized $\sin \theta$ theorem~\cite{Wedin}. To state such a result, it is instrumental to introduce an appropriate rewriting of the SVD's in \eqref{Ufac} and \eqref{svd_perturbed}. By exploiting that matrices $\Xi$, $K$, and $\tilde \Xi$, $\tilde K$ are unitary, it follows that
\begin{equation*}
U = U_1 + U_0,\quad \tilde U = \tilde U_1 + \tilde U_0,
\end{equation*}
where $U_s = \Xi_s \Lambda_s K_s^T$, $\tilde U_s = \tilde \Xi_s \tilde \Lambda_s \tilde K_s^T$, for $s=0, 1$,
with
$K_1 = [\mathbf{k}_1, \ldots, \mathbf{k}_r]$,
$K_0 = [\mathbf{k}_{r+1}, \ldots, \mathbf{k}_{mp}]$,
$\Xi_1 = [\boldsymbol{\xi}_1, \ldots, \boldsymbol{\xi}_r]$,  
$\Xi_0 = [\boldsymbol{\xi}_{r+1},\ldots, \boldsymbol{\xi}_{N_h}]$,
\begin{equation*}
\Lambda_1=\left\{
\begin{array}{ll}
\begin{bmatrix}{\rm diag}(\lambda_1, \ldots, \lambda_r, 0, \ldots,0) \\ \mathbf{0}_{(N_h -r)\times N_h}\end{bmatrix} \qquad &\text{if} \ N_h \ge mp \\[5mm]
\begin{bmatrix}{\rm diag}(\lambda_1, \ldots, \lambda_r, 0, \ldots,0), \mathbf{0}_{mp\times (mp -r)}\end{bmatrix} & \text{if} \ N_h < mp,\\[1mm]
\end{array}
\right.
\end{equation*}
$$\Lambda_0=\left\{
\begin{array}{ll}
\begin{bmatrix}{\rm diag}(0, \ldots, 0,\lambda_{r+1}, \ldots, \lambda_{\min\{N_h, mp\}})\\\mathbf{0}_{(N_h -r)\times N_h}\end{bmatrix} \quad &\text{if} \ N_h \ge mp \\[5mm]
\begin{bmatrix}{\rm diag}(0, \ldots, 0,\lambda_{r+1}, \ldots, \lambda_{\min\{N_h, mp\}}),\mathbf{0}_{mp\times (mp -r)}\end{bmatrix} &\text{if} \ N_h < mp,\\[1mm]
\end{array}
\right.
$$
matrices $\tilde K_1$, $\tilde K_0$, $\tilde \Xi_1$, $\tilde \Xi_0$,
$\tilde \Lambda_1$, $\tilde \Lambda_0$ being defined accordingly.\\
Moreover, we denote by ${\mathcal C}_1$ and ${\mathcal C}_2$ two generic Euclidean subspaces of $\mathbb{R}^{N_h}$, and by $P_{{\mathcal C}_1}$, $P_{{\mathcal C}_2}\in \mathbb{R}^{N_h}$ the associated orthogonal projection operators. Thus, the angle $\theta$ between a vector $\mathbf{x}\in \mathbb{R}^{N_h}$ and the subspace ${\mathcal C}_1$ (which, among all the mathematically equivalent formulations, by convention is always taken acute and positive) can be defined by
\begin{equation}
    \sin\theta(\mathbf{x},{\mathcal C}_1)=\min_{\mathbf{y}\in {\mathcal C}_1} \lVert \mathbf{x}-\mathbf{y} \rVert_2,
\end{equation}
with $\lVert\cdot\rVert_2$ denoting the Euclidean norm of a vector, and with $\lVert \mathbf{x}\rVert_2=1$.
From the projection theorem, it follows that 
\begin{equation}
    \min_{\mathbf{y}\in {\mathcal C}_1} \lVert \mathbf{x}-\mathbf{y} \rVert_2 = \lVert (I-P_{{\mathcal C}_1}\mathbf{x}) \rVert_2.
\end{equation}
When considering the angle between the subspaces ${\mathcal C}_1$ and ${\mathcal C}_2$, it is common to define
\begin{equation}
    \sin \theta ({\mathcal C}_1, {\mathcal C}_2) = \lVert (I-P_{{\mathcal C}_2})P_{{\mathcal C}_1} \rVert_S.
\end{equation}

Now, 
if there exists a pair of values $\gamma>0$, $\delta>0$ such that 
\begin{equation}
    \lambda_{\min}(\tilde U_1) \ge \gamma + \delta, \quad 
    \lambda_{\max}(U_0) \le \gamma,
\end{equation} 
with $\lambda_{\min}(\tilde U_1)$ and $\lambda_{\max}(U_0)$ the minimum and the maximum eigenvalue of $\tilde U_1$ and $U_0$, respectively, the generalized $\sin \theta$ theorem~\cite{Wedin} ensures that the following perturbation bound on the column space of $\tilde \Xi$ holds
\begin{equation}\label{bound1}
    \sin \theta \big (R(\tilde \Xi_1),R( \Xi_1) \big)  \le \frac{\lVert E \rVert_S}{\delta}={\mathcal P}_{B1},
\end{equation}
with $R(W)$ the column space associated with the generic matrix $W$. 
The relation \eqref{bound1} is not always computationally convenient to provide meaningful insight on the actual perturbation triggered by the noise in the data. Indeed, computing the perturbation between the two column spaces $R(\tilde \Xi_1)$ and $R( \Xi_1)$ 
might be unpractical from a computational view point, due to the infinite vectors to be spanned.\\
In order to make the control in \eqref{bound1} more  computationally convenient, i.e., to restrict the dimensionality of the spaces to be spanned, we consider the alternative inequality 
\begin{equation}\label{sin1}
     \sin \theta \Big(R(\tilde \Xi_1),R( \Xi_1) \Big)  \le  \sin \theta (\tilde{\boldsymbol{\xi}}_1,\boldsymbol{\xi}_1 ) ={\mathcal P}_{B2},
\end{equation}
with $\boldsymbol{\xi}_1$ and $\tilde{\boldsymbol{\xi}}_1$ denoting the first left singular vectors of the response matrix $U$ and $\tilde U$, respectively. 
The result in \eqref{sin1} can be easily proved, moving from the following chain of inequalities: 
\begin{align*}
 \sin\theta (\tilde{\boldsymbol{\xi}}_1,\boldsymbol{\xi}_1 )  & \ge  \sin \theta (\tilde{\boldsymbol{\xi}}_1,\text{span}\{\boldsymbol{\xi}_1, \boldsymbol{\xi}_2\} )  \ge \sin \theta (\tilde{\boldsymbol{\xi}}_1,\text{span}\{\boldsymbol{\xi}_1, \boldsymbol{\xi}_2, \boldsymbol{\xi}_3\} ) \\[1mm]
& \ge \ldots
\ge \sin \theta (\tilde{\boldsymbol{\xi}}_1,R(\Xi) ) 
\ge  \sin \theta (\text{span}\{\tilde{\boldsymbol{\xi}}_1, \tilde{\boldsymbol{\xi}}_2\},R(\Xi) ) \\[1mm]
& \ge \sin \theta (\text{span}\{\tilde{\boldsymbol{\xi}}_1, \tilde{\boldsymbol{\xi}}_2, \tilde{\boldsymbol{\xi}}_3\},R(\Xi) ) 
\ge \ldots \\[1mm]
& 
\ge \sin \theta \Big(R(\tilde \Xi_1),R( \Xi_1) \Big) .
\end{align*}

Bounds \eqref{weyl}, \eqref{mirsky}, \eqref{bound1}, together with the new one in \eqref{sin1}, will be used in the next section 
to assess the sensitivity 
both of the interpolation-HiPOD and of the regression-HiPOD
to the noise level. 

\section{Numerical results}\label{num_sec}
In this section we compare the interpolation-HiPOD with the regression-HiPOD methods. 
The comparison in carried out on two test cases, in terms of accuracy and robustness to the noise level. For the former method, we resort to the PCH interpolant, while the latter approach is assessed both with the polynomial and with the Gaussian process regression.

\subsection{Test case 1}\label{tc1}
We choose as reference setting the ADR problem \eqref{parametrized_referenc}-\eqref{forms} solved on the rectangular domain $\Omega=(0,6)\times(0,1)$, after setting the problem data to
\begin{equation}\label{setting1}
\begin{array}{c}
\mu(x, y)=0.24, \quad \mathbf{b}(x, y)=[5,  \sin(6x)]^T, \quad \sigma(x, y) = 0.1,\\[2mm]
f(x,y)=10\chi_{E_1}(x, y)+10\chi_{E_2}(x, y), 
\end{array}
\end{equation}
with $\chi_{\omega}$ the characteristic function associated with the generic region $\omega\subset \mathbb R^2$, 
$E_1$ and $E_2$ the ellipsoidal areas in $\Omega$ given by $\{ (x,y):(x-0.75)^2+0.4(y-0.25)^2 < 0.01\}$  and 
$\{ (x,y):(x-0.75)^2+0.4(y-0.75)^2 <0.01 \}$, respectively. 
The ADR problem is completed with a homogeneous Neumann data on
$\Gamma_N=\{ (x,y): x=6, 0\le y \le 1\}$, while a homogeneous Dirichlet condition is assigned on $\Gamma_D=\partial \Omega \setminus \Gamma_N$, so that space $V$ in \eqref{parametrized_referenc} coincides with 
$H^1_{\Gamma_D}(\Omega)$.\\
From a modeling viewpoint, this setting can be adopted to simulate the propagation of a pollutant released by two localized sources within a straight channel, under the effect of a sinusoidal horizontal convection.

\begin{figure}
\includegraphics[width=1.002\textwidth]{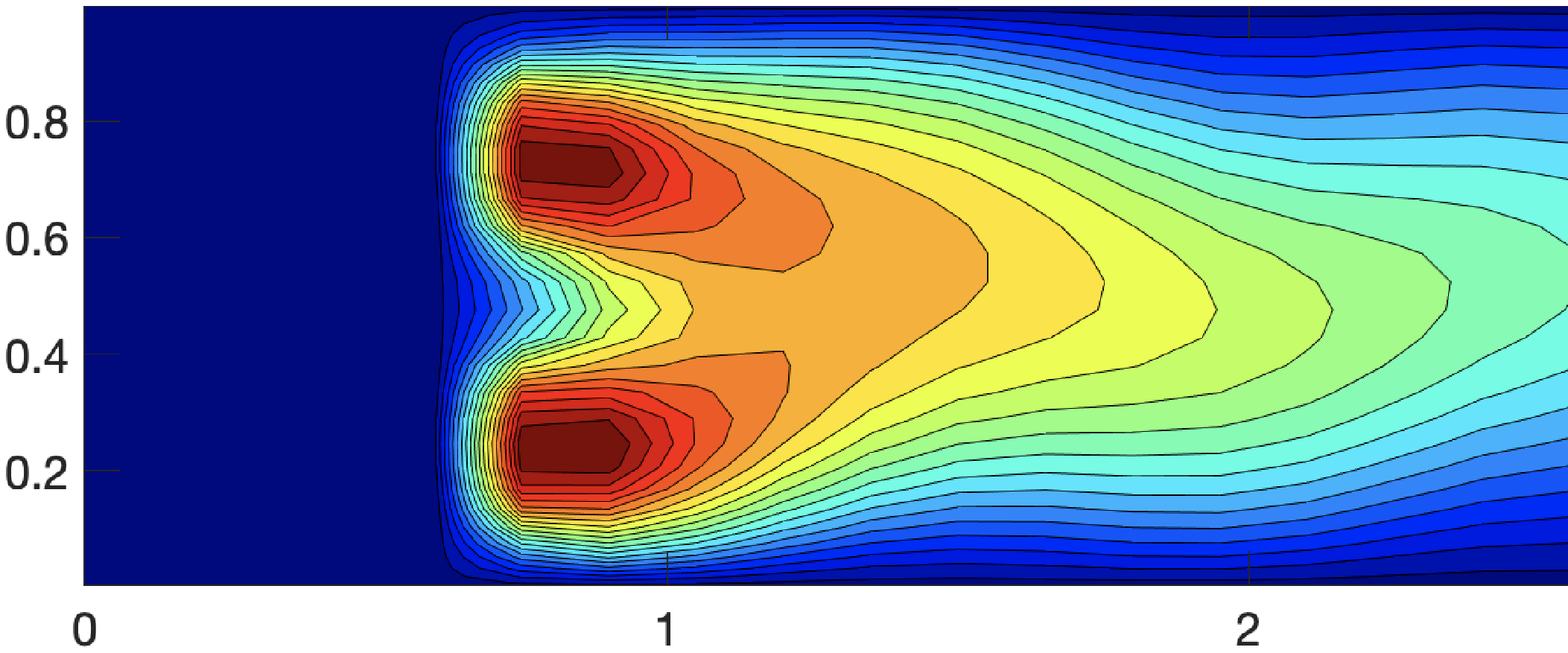}
\includegraphics[width=0.998\textwidth]{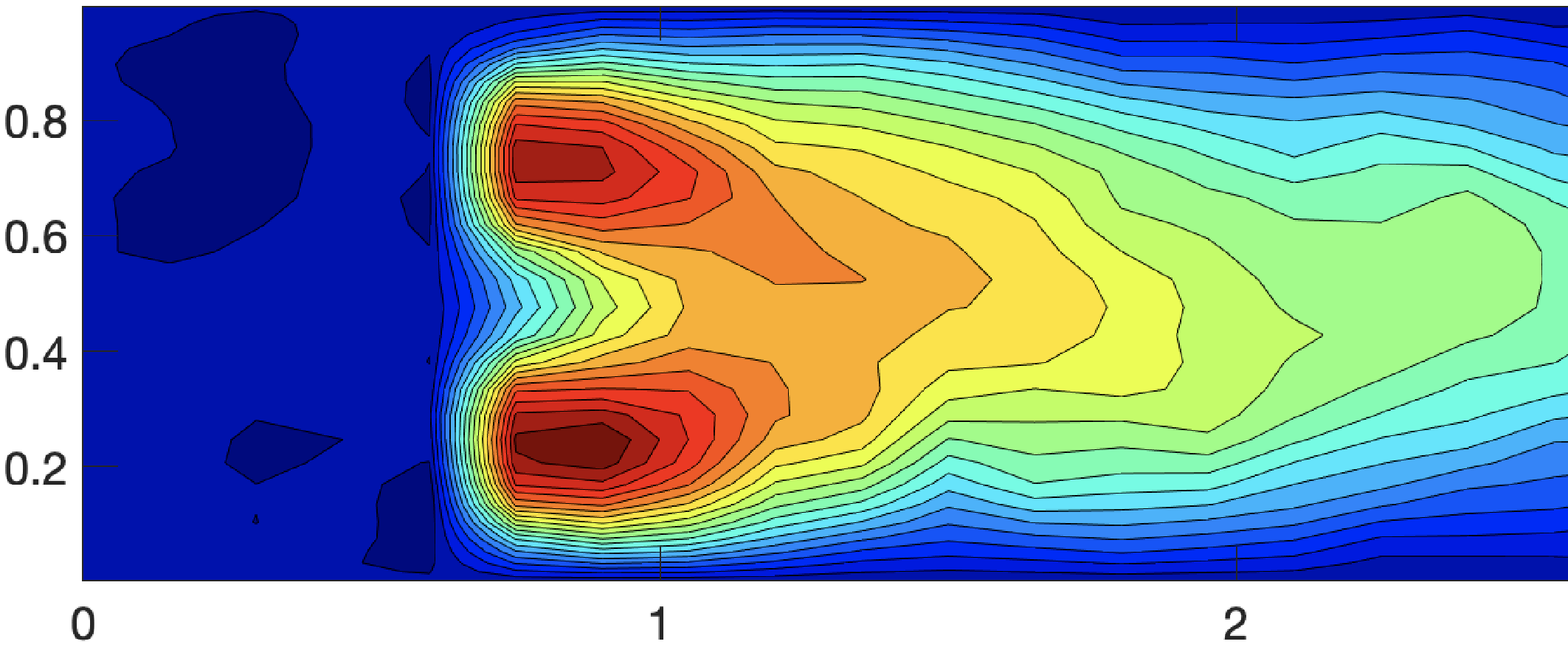}
\includegraphics[width=\textwidth]{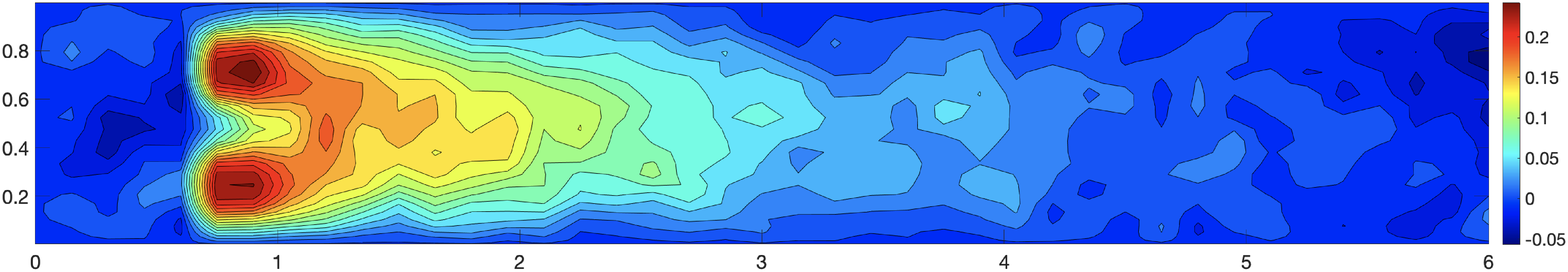}
\includegraphics[width=0.997\textwidth]{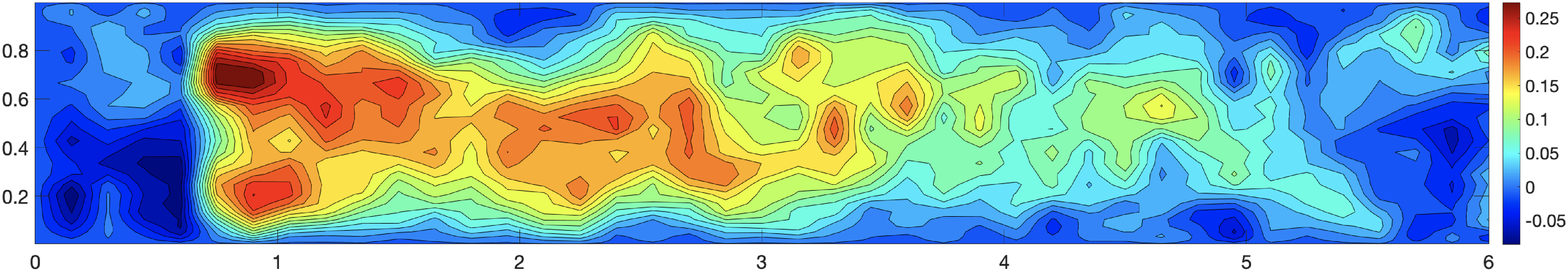}
\caption{Test case 1. HiMod reference solution (top panel);  interpolation-HiPOD approximation for different noise levels: $\eta = 0.01$, $0.1$, $0.25$ (second-fourth row).}
\label{test1_interp}
\end{figure}

We parametrize the ADR problem with respect to the diffusivity $\mu$.
The offline phase is set up so that the coefficient $\mu$ uniformly spans the range $[0.2, 0.8]$. All the other problem data remain the same as in \eqref{setting1} through the entire offline phase.
We hierarchically reduce $100$
ADR problems, after discretizing the main dynamics with linear finite elements on a uniform partition of $\Omega_{1D}$ into $40$ sub-intervals, while using $20$ sinusoidal modal basis functions to approximate  
the transverse dynamics.             
To investigate the robustness to the noise of the interpolation- and of the regression-HiPOD procedures, we carry out different offline phases, where the 
source term $f$ is 
injected by different levels of white noise, thus changing $f$ into $f+\eta$,
with $\eta = 0.01$, $0.05$, $0.1$, $0.25$. Independently of the selected noise, the online phase is used to recover the reference setting in \eqref{setting1}, i.e., to reconstruct the HiMod solution for $\mu^* = 0.24$. The top panel in Fig.~\ref{test1_interp} shows the approximation yielded by the HiMod discretization adopted for the offline phase. 
The resolution of such a discretization is sufficiently accurate to capture the oscillatory dynamics induced by the sinusoidal field, 
together with the presence of the two localized sources in $E_1$ and $E_2$. Such a discretization represents the reference solution HiPOD approximations will be compared with.

A unique threshold tolerance, $\varepsilon$, is used
to automatically select the number of POD modes to be preserved at the first and at the second stage of the directional method, which is equivalent to setting $\epsilon_1=\epsilon_2=\epsilon$ in \eqref{tolerances} (we refer the reader to~\cite{lupopasini2022} for a thorough investigation about the interplay between tolerances $\epsilon_1$ and $\epsilon_2$).  Table~\ref{test_1_singular_values} collects the distribution of the 
number $\mu_j$ of the left singular vectors retained at each finite element node, for different choices of the tolerance $\varepsilon$ (by rows) and of the noise $\eta$ (by columns).
As expected, number $\mu_j$ increases for larger and larger values both of  $\varepsilon$ and $\eta$. 
Indeed, in the former case, the procedure is intrinsically requested to retain more information about the variability of the offline dataset. In the latter case, the quality of information retained by the offline solution deteriorates with increasing values of noise. This leads the HiPOD procedure to retain more singular vectors with respect to the case of a low (or of the absence of) noise, in order to capture the same amount of information from the offline data. 
\begin{table}\sffamily
\begin{tabular}{c|c|c}
\toprule
$\eta=0.01$ & $\eta=0.1$ & $\eta=0.25$ \\ 
\midrule \hspace*{-.3cm}
 \includegraphics[
 width=0.31\textwidth]{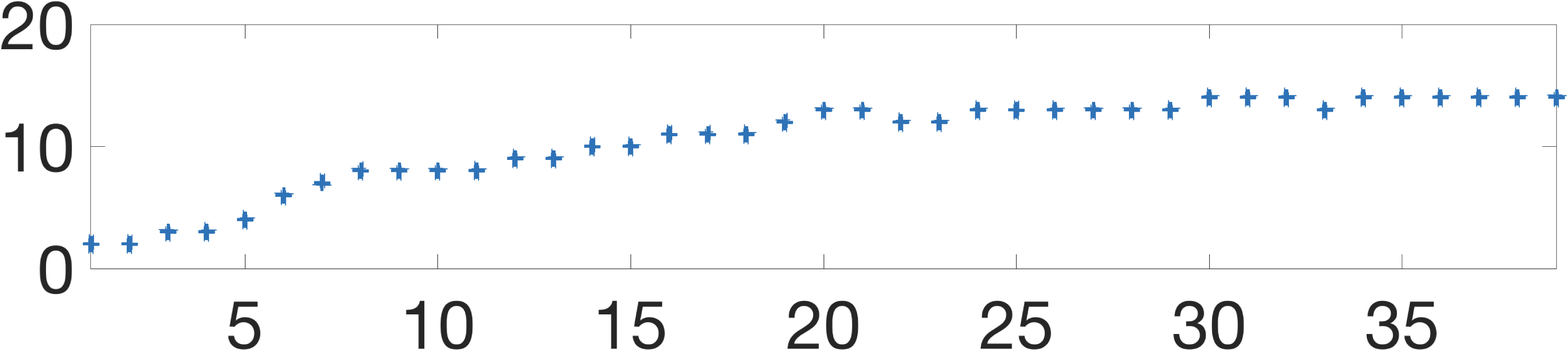} & \hspace*{-.2cm} \includegraphics[
 width=0.31\textwidth]{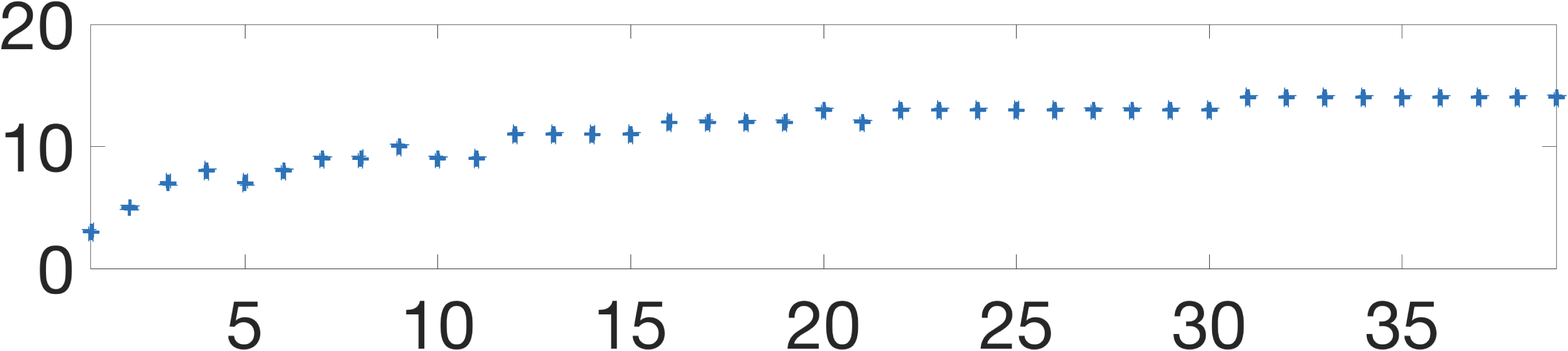} & \hspace*{-.2cm} \includegraphics[
 width=0.31\textwidth]{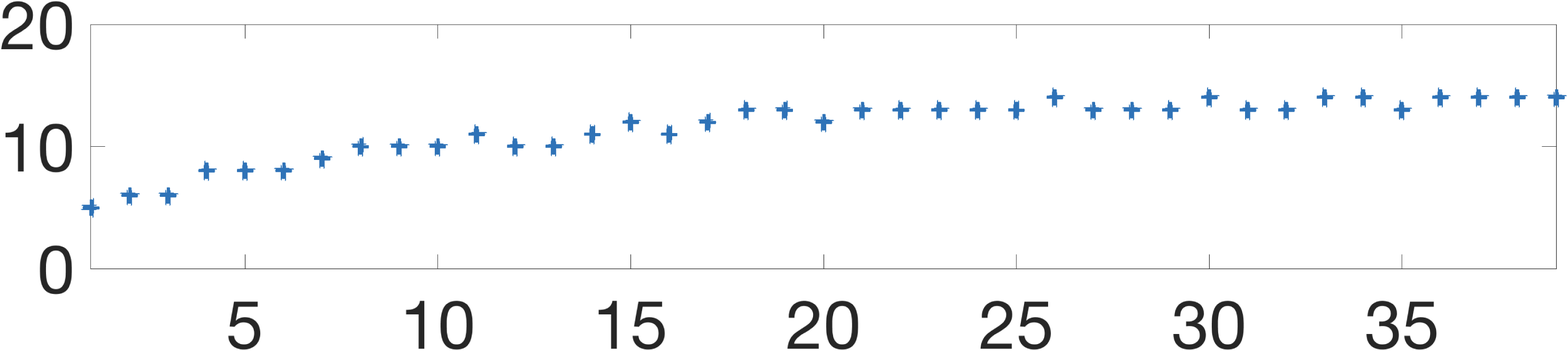} \\
 & & \\
 \hspace*{-.3cm} \includegraphics[width=0.31\textwidth]{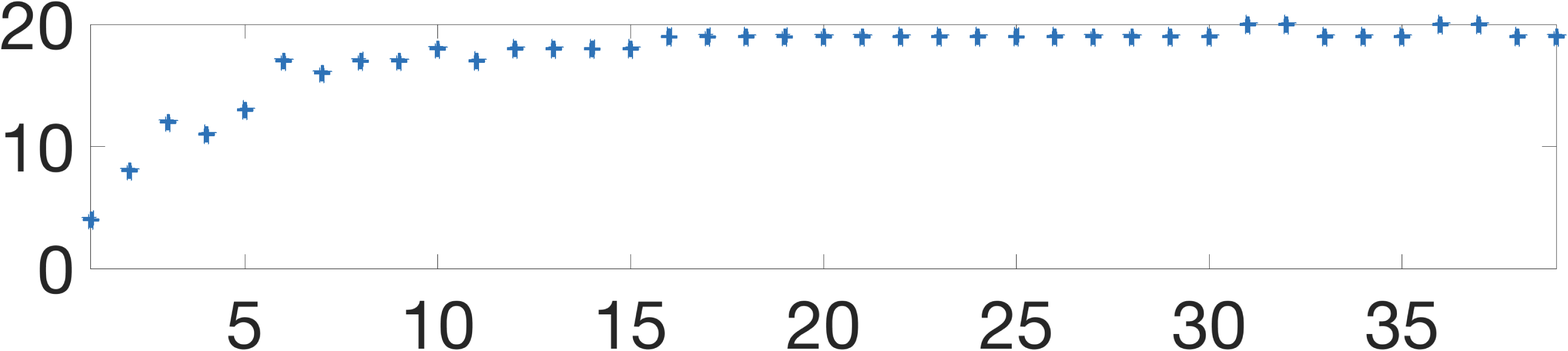} & \hspace*{-.2cm} \includegraphics[width=0.31\textwidth]{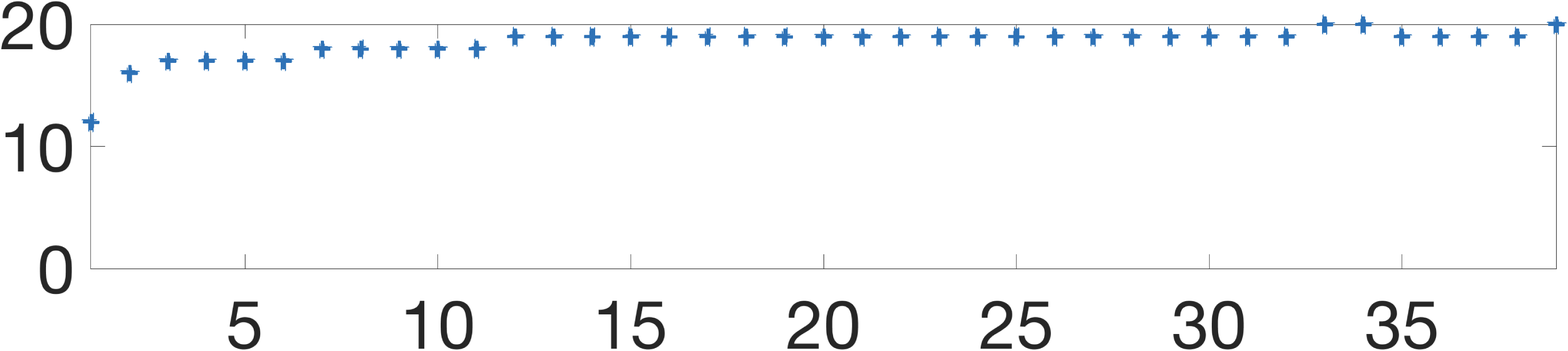} & \hspace*{-.2cm} \includegraphics[width=0.31\textwidth]{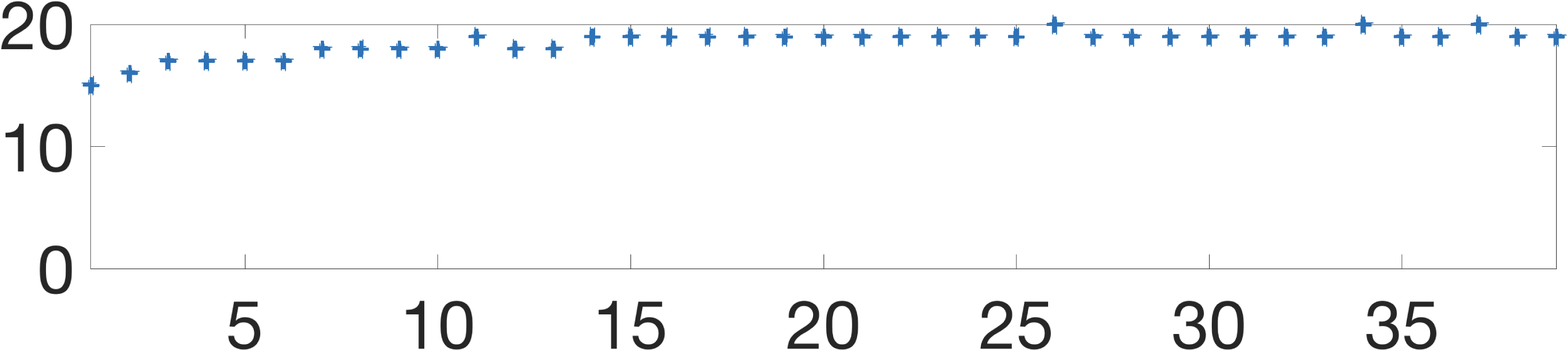} \\ 
 &&\\
 \hspace*{-.3cm} \includegraphics[width=0.31\textwidth]{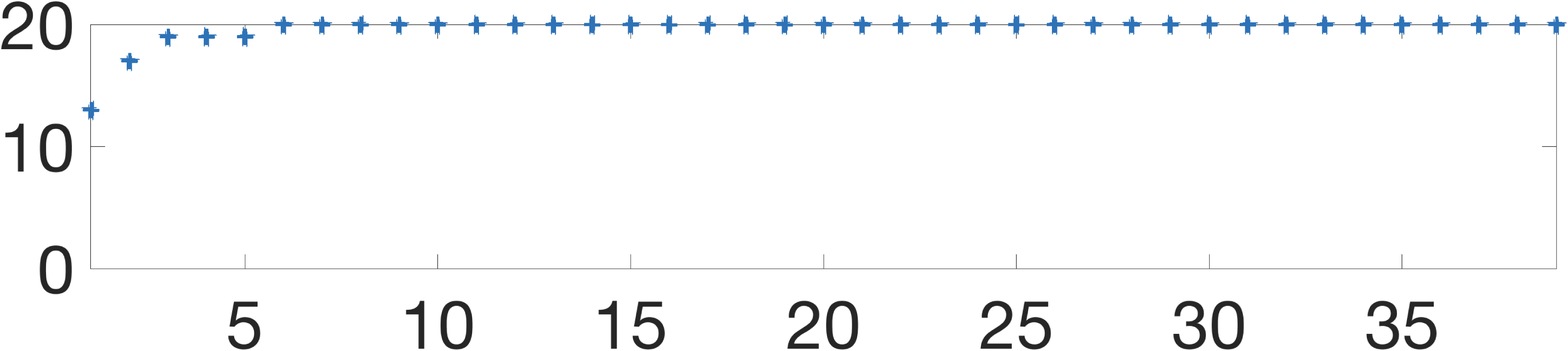} & \hspace*{-.2cm} \includegraphics[width=0.31\textwidth]{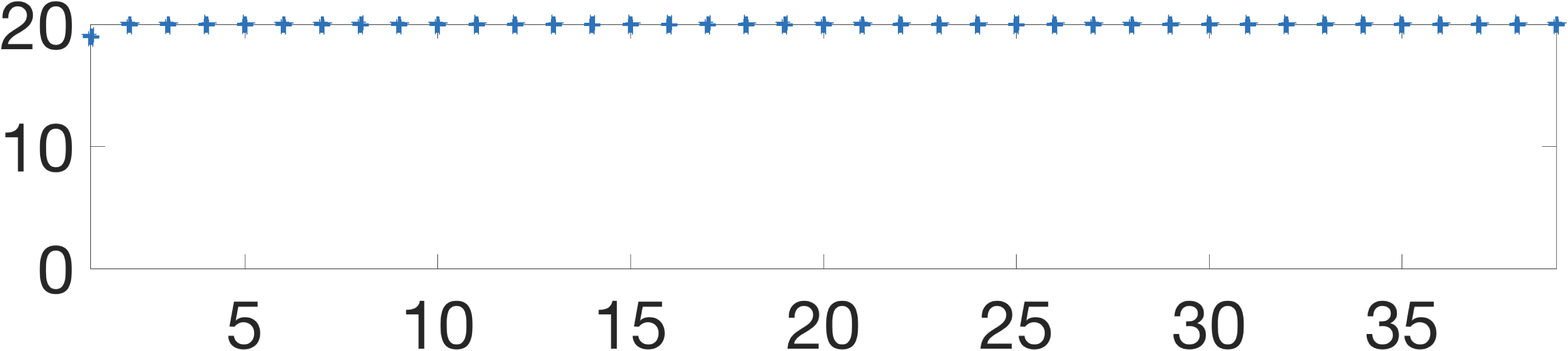} & \hspace*{-.2cm} \includegraphics[width=0.31\textwidth]{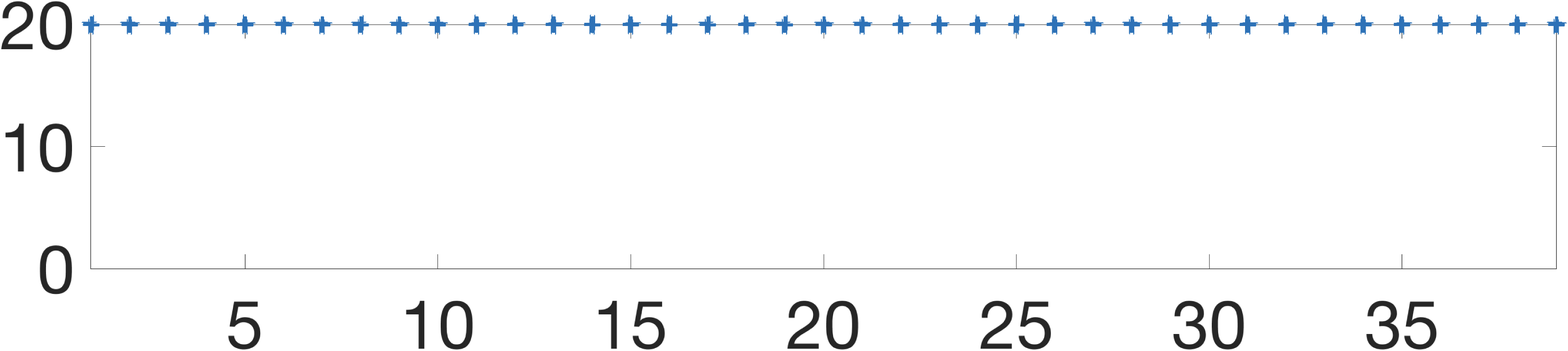} \\ 
\bottomrule 
\end{tabular}
\caption{Test case 1. Distribution of the number of the POD modes retained at each finite element node for $\varepsilon=0.9$ (first row), $0.99$ (second row), $0.999$ (third row), and for a different noise $\eta$.
The supporting fiber is discretized by $41$ uniformly distributed nodes.}
\label{test_1_singular_values}
\end{table}

Figure~\ref{test1_interp} 
displays the interpolation-HiPOD approximation reconstructed from offline data with an increasing level of noise, i.e., $\eta = 0.01$, $\eta =0.1$, $\eta = 0.25$
(the contourplot for $\eta =0.05$ is omitted since it is very similar to the one associated with $\eta=0.1$). As expected, the quality of the HiPOD solution deteriorates with increasing values of $\eta$. In particular, the wake behind the sources is completely lost by the HiPOD approximation for $\eta = 0.25$.

\begin{figure}
\includegraphics[width=1.005\textwidth]{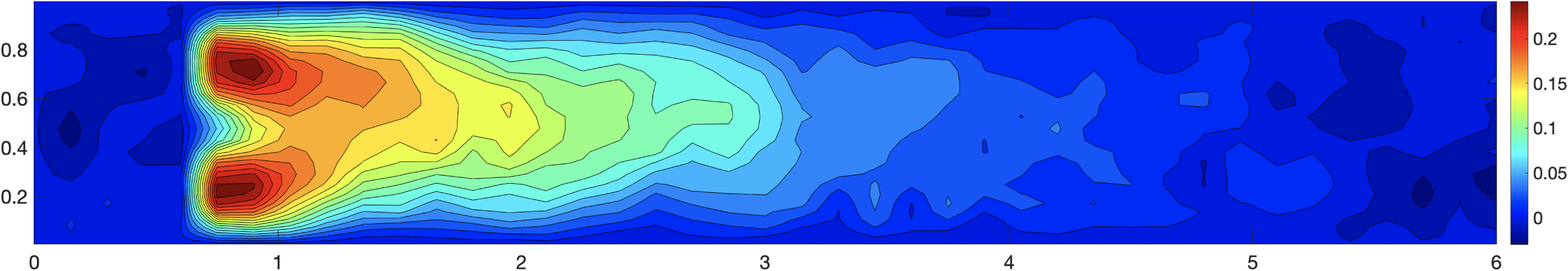}
\includegraphics[width=\textwidth]{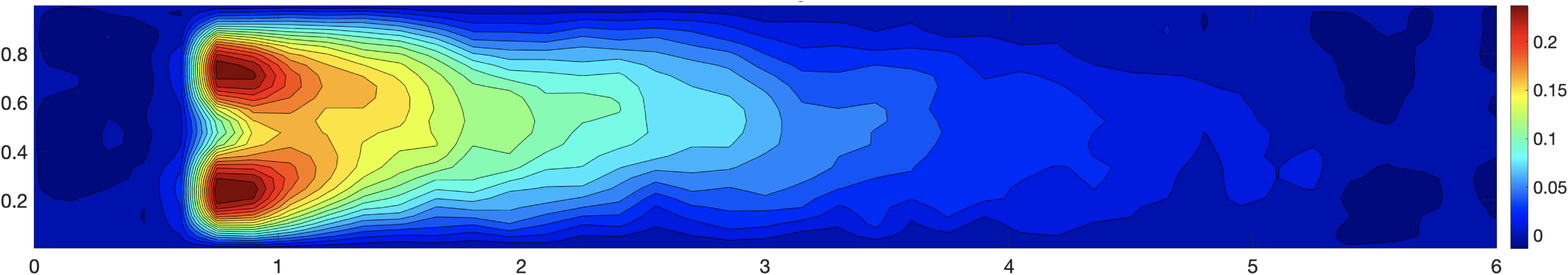}
\caption{Test case 1. Regression-HiPOD approximation for a cubic polynomial fitting (top) and for a Gaussian process regression (bottom); the noise level is set to $\eta = 0.25$.}\label{test1_fitting}
\end{figure}

Figure~\ref{test1_fitting}
highlights the benefits obtained by replacing the interpolation step in the HiPOD online phase with a regressione process. In particular, the two panels 
show the regression-HiPOD approximation when resorting to a cubic polynomial fitting and to a Gaussian process regression, for the largest level of noise analyzed in Fig.~\ref{test1_interp} (i.e., $\eta=0.25$). Both fitting models clearly outperform the  interpolation-HiPOD. In particular, the Gaussian process generates a solution of a better quality with respect to the one yielded by the cubic polynomial fitting. 
\begin{table}[h]
\begin{center}
\begin{tabular}{|c|c|c|c|c|c|}
\hline
  & & $\eta = 0.01$ & $\eta = 0.05$ & $\eta = 0.1$ & $\eta = 0.25$\\
  \hline
  PCH & $L^2(\Omega)$-norm & 0.0222& 0.0960& 0.2258 & 0.5521 
 \\[1mm] \cline{2-6}
   & $H^1(\Omega)$-norm & 0.0434& 0.2209& 0.4460 & 1.1431 
 \\[1mm]
  \hline
   CP & $L^2(\Omega)$-norm & 0.0060 & 0.0247 & 0.0648 & 0.1414 
 \\[1mm] \cline{2-6}
   & $H^1(\Omega)$-norm & 0.0132 & 0.0592 & 0.1282 & 0.3913 
 \\[1mm]
  \hline
  GP & $L^2(\Omega)$-norm & 0.0096& 0.0229 & 0.0742 & 0.1521 
 \\[1mm] \cline{2-6}
   & $H^1(\Omega)$-norm & 0.0169& 0.0559 & 0.1087 & 0.2730 
 \\[1mm]
  \hline
\end{tabular}
\end{center}
\caption{Test case 1. Relative modeling error associated with the interpolation-HiPOD (PCH) and with the cubic polynomial (CP) and the Gaussian process (GP) regression-HiPOD, for different levels of noise.}
\label{tab_test1_hipod}
\end{table}

This trend is confirmed by the values in Table~\ref{tab_test1_hipod}, which gathers the $L^2(\Omega)$- and the $H^1(\Omega)$-norm of the relative modelig error obtained 
when replacing the reference HiMod discretization in Fig.~\ref{test1_interp} (top) with the
interpolation- rather than the regression-HiPOD approximation. To ensure a high accuracy to the HiPOD solutions, we have set the
thresholds driving the selection of the POD bases in \eqref{tolerances} very close to $1$, picking $\epsilon_1=\epsilon_2=\epsilon=0.9999$. The interpolation-HiPOD scheme is compared with the regression-HiPOD approach when both a cubic polynomial fitting and a Gaussian process are used, for the four levels of noise considered above. 
Although the general trend confirms that the error increases with the noise, cubic polynomial fitting and Gaussian process regression outperform considerably the interpolation-HiPOD, in particular for high noise levels. 

\begin{table}[h]
\begin{center}
\begin{tabular}{|c|c|c|c|c|}
\hline
 & $\eta = 0.01$ & $\eta = 0.05$ & $\eta = 0.1$ & $\eta = 0.25$\\
 \hline
  ${\mathcal P}_{B1}$ & 2.2819 & 3.3340 & 2.6980 & 3.9770
 \\ 
  \hline
   ${\mathcal P}_{B2}$ & 0.0032 & 0.0072 & 0.0173 & 0.0517
 \\   
 \hline
\end{tabular}
\end{center}
\caption{Test case 1. Perturbation bounds for the noise propagation in the response matrix.}
\label{tab_test1_sin}
\end{table}

Finally, we quantify the  perturbation bounds in \eqref{bound1} and \eqref{sin1}, for the considered four levels of noise. Actually,   Table~\ref{tab_test1_sin}
shows that the new bound allows us to gain up to three orders of magnitude with respect to quantifier ${\mathcal P}_{B1}$, thus 
offering an effective tool to evaluate the noise propagation on the response matrix $U$ in \eqref{response_matrix}.

\subsection{Test case 2}\label{sec2}
As a second reference setting, we hierarchically reduce the ADR problem in \eqref{parametrized_referenc}-\eqref{forms} identified by the data 
\begin{equation}\label{setting2}
\begin{array}{c}
\mu(x, y)=0.24, \quad \mathbf{b}(x, y)=[20,  2  \sin(6x)]^T, \quad \sigma(x, y) = 0.1,\\[2mm]
f(x,y)=1000\chi_{R_1}(x, y) + 1000\chi_{R_2}(x, y), 
\end{array}
\end{equation}
on the same domain, $\Omega=(0,6)\times (0,1)$, as in Test case 1, with $R_1=\{ (x,y): 1<x<2, \; 0<y<0.1 \}$ and $R_2=\{ (x,y): 1<x<2, \; 0.9<y<1.0 \}$ two rectangular regions of interest.
The problem is completed by homogeneous Neumann data on 
$\Gamma_N=\{ (x,y): x=6, 0\le y \le 1\}$ and by a homogeneous Dirichlet condition on $\Gamma_D=\partial \Omega \setminus \Gamma_N$, so that we have $V\equiv H^1_{\Gamma_D}(\Omega)$ in \eqref{parametrized_referenc}.\\
We can adopt this configuration to model, for instance,
the transport of a drug released by a medical stent applied to the walls of a cardiovascular vessel, under the effect of incompressible fluid convection.

\begin{figure}
\includegraphics[width=1.002\textwidth]{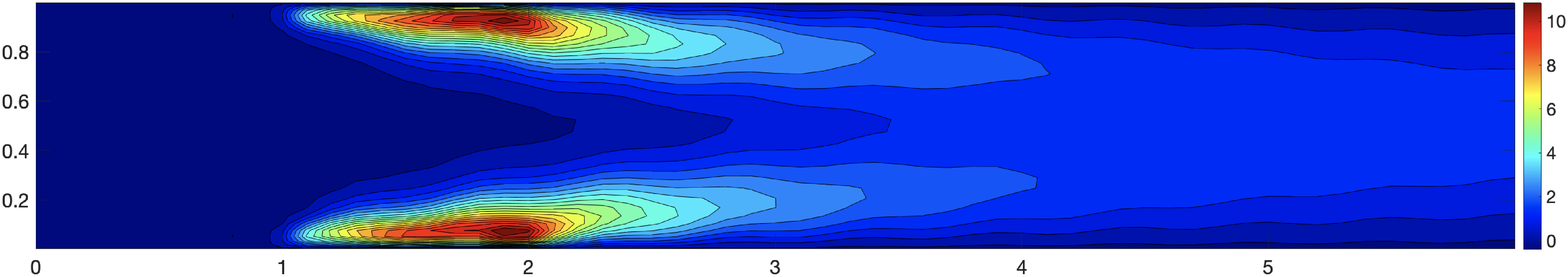}
\includegraphics[width=\textwidth]{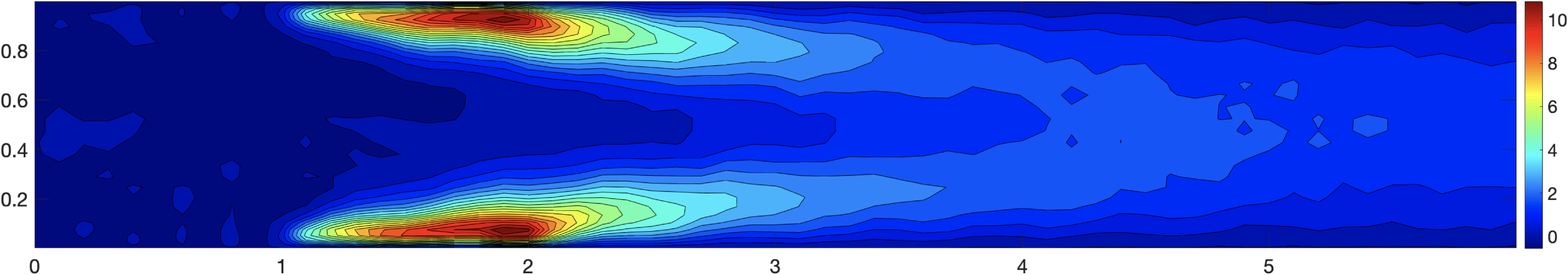}
\includegraphics[width=\textwidth]{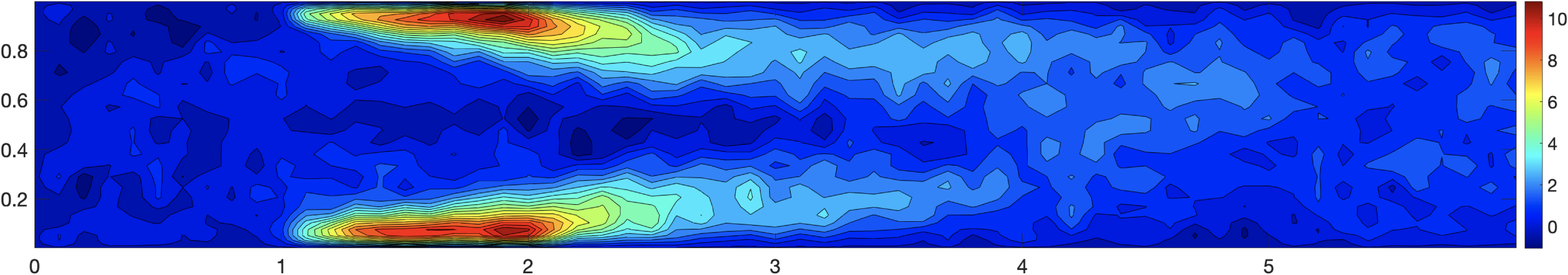}
\includegraphics[width=\textwidth]{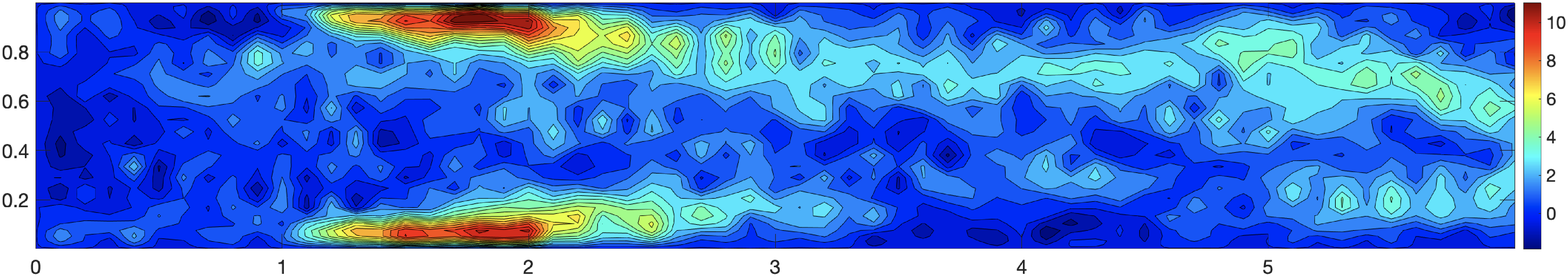}
\includegraphics[width=\textwidth]{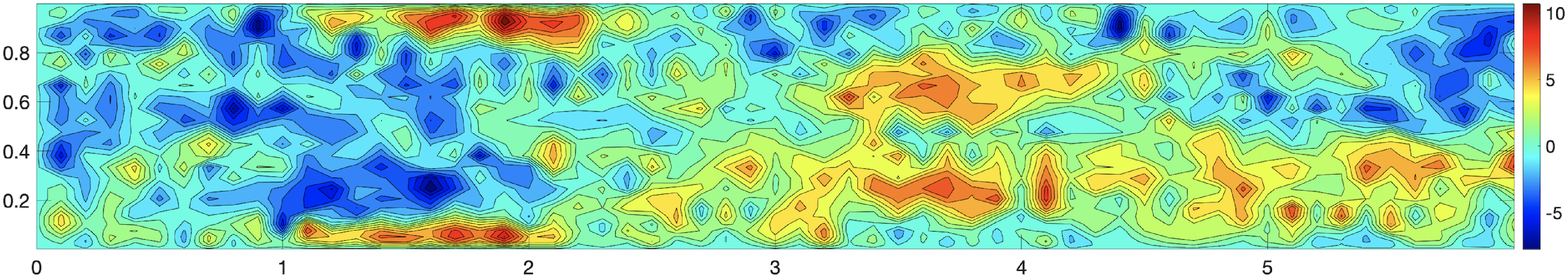}
\caption{Test case 2. HiMod reference solution (top panel); interpolation-HiPOD approximation for different noise levels: $\eta = 0.01$, $0.05$, $0.1$, $0.25$ (second-fifth row). }\label{test2_interp}
\end{figure}

The HiMod discretization adopted to build the reference solution as well as to perform the offline phase of the HiPOD procedure employs linear finite elements along $\Omega_{1D}$ after subdiving the supporting fiber into $60$ uniform subintervals, and a modal basis consisting of $20$ sinusoidal functions to 
capture the dynamics along the transverse direction. The top panel in Fig.~\ref{test2_interp} displays the contour plot of the HiMod reference solution.
We clearly distinguish the drug release in the regions $R_1$ and $R_2$, together with the  transport of the medicine along the pipe.
The offline phase is driven by the diffusivity, which is identified with parameter $\alpha$. In particular, we uniformly cover the range $[0.2, 0.8]$ with $100$ samples, while keeping the same values as in \eqref{setting2} for all the other problem data.

The performance of the interpolation- and of the regression-HiPOD approaches in the presence of noise is analyzed by affecting the forcing term in \eqref{forms} with an increasing white noise $\eta$, set to $\eta=0.01$, $0.05$, $0.1$, $0.25$, respectively. 
For each choice of $\eta$, we replicate the offline phase set above, before predicting the HiMod approximation associated with parameter $\mu^*=0.24$ (i.e., with the reference configuration) in the online phase.

The selection of the POD bases $\{ {\boldsymbol \xi}_j\}_{j=1}^L$ and $\{ {\bf r}_j^k\}_{k=1}^{\mu_j}$, for $j=1, \ldots, L$, is driven by the variance-based criteria in \eqref{tolerances}, for a unique choice of the threshold tolerance (i.e., for $\epsilon_1=\epsilon_2=\epsilon$). 
The trend of the number $\mu_j$ of the left singular vectors selected by the HiPOD approach
at the finite element nodes for different values of $\epsilon$ and $\eta$ is very similar to the one in Table~\ref{test_1_singular_values} (and, consequently, skipped for shortness).

\begin{figure}
\includegraphics[width=1.005\textwidth]{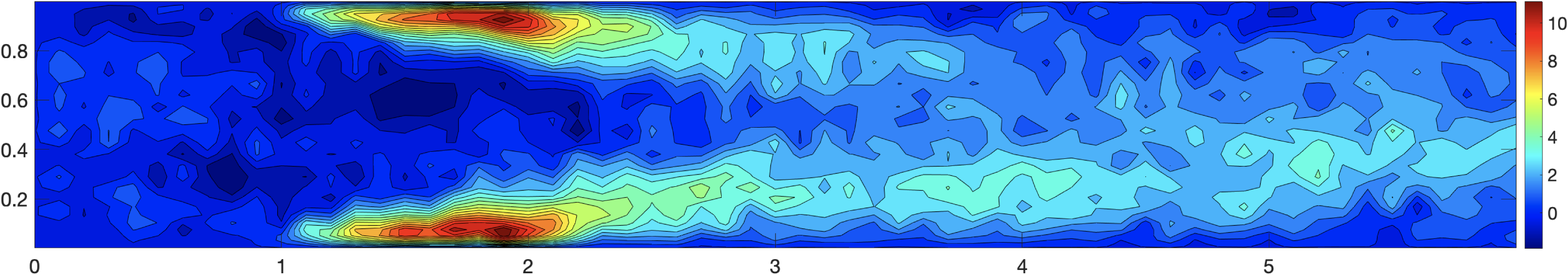}
\includegraphics[width=\textwidth]{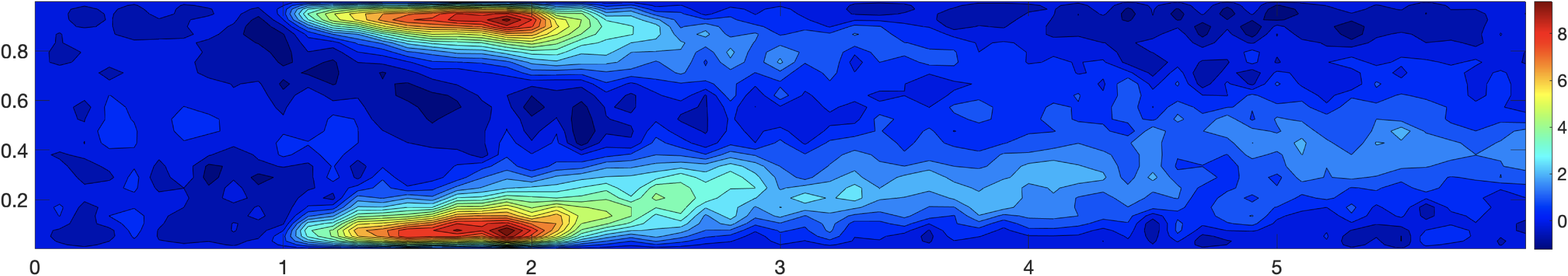}
\caption{Test case 2. Regression-HiPOD approximation for a cubic polynomial fitting (top) and for a Gaussian process regression (bottom); the noise level is set to $\eta = 0.25$.}\label{test2_fitting}
\end{figure}

Figure~\ref{test2_interp} shows the interpolation-HiPOD solution for the four noise levels $\eta=0.01$, $0.05$, $0.1$, $0.25$. 
As for Test case 1, the quality of the HiPOD approximation deteriorates very quickly when increasing the level of noise.
In particular, for $\eta=0.25$, the HiPOD solution is fully noisy, the problem dynamics being completely lost.\\
On the contrary, the regression-HiPOD discretization significantly outperforms the approximation quality provided by the interpolation.
Figure~\ref{test2_fitting} shows such an improvement when resorting to a cubic polynomial fitting and to a Gaussian process regression. Analogously to Fig.~\ref{test1_fitting},
Gaussian process regression yields a solution of a better quality than cubic polynomial fitting.

These qualitative considerations are confirmed by the values in Table~\ref{tab_test2_hipod}, which gathers the $L^2(\Omega)$- and the $H^1(\Omega)$-norm of the relative modeling error associated with both the interpolation- and the regression-HiPOD approximations, for the considered noise levels.
The results corroborate what already remarked for the first test case, namely cubic polynomial fitting and Gaussian process regression outperform the interpolation-HiPOD for high noise levels.
\begin{table}[h]
\begin{center}
\begin{tabular}{|c|c|c|c|c|c|}
\hline
  & & $\eta = 0.01$ & $\eta = 0.05$ & $\eta = 0.1$ & $\eta = 0.25$\\[1mm]
  \hline
  PCH & $L^2(\Omega)$-norm & 0.0446& 0.1917& 0.4310 & 1.0589 
 \\[1mm] \cline{2-6}
   & $H^1(\Omega)$-norm & 0.0626& 0.3177& 0.6723 & 1.6280 
 \\[1mm]
  \hline
   CP & $L^2(\Omega)$-norm & 0.0163 & 0.0887 & 0.1408 & 0.3067 
 \\[1mm] \cline{2-6}
   & $H^1(\Omega)$-norm & 0.0213 & 0.0928 & 0.1875 & 0.4052
 \\[1mm]
  \hline
  GP & $L^2(\Omega)$-norm & 0.0137& 0.0832 & 0.1617 & 0.3593
 \\[1mm] \cline{2-6}
   & $H^1(\Omega)$-norm & 0.0182& 0.0982& 0.1758 & 0.3703 
 \\[1mm]
  \hline
\end{tabular}
\end{center}
\caption{Test case 2.
Relative modeling error associated with the interpolation-HiPOD (PCH) and with the cubic polynomial (CP) and the Gaussian process (GP) regression-HiPOD, for different levels of noise.}
\label{tab_test2_hipod}
\end{table}

As a last check, we compute the quantifiers ${\mathcal P}_{B1}$ and ${\mathcal P}_{B2}$ of the noise propagation in the response matrix defined in \eqref{bound1} and \eqref{sin1}, respectively (see Table~\ref{tab_test2_sin}).
Also for this test case, the bound in~\eqref{bound1} is
relatively large, and the discrepancy between ${\mathcal P}_{B1}$ and ${\mathcal P}_{B2}$ grows with the level of noise.
\begin{table}[h]
\begin{center}
\begin{tabular}{|c|c|c|c|c|}
\hline
 & $\eta = 0.01$ & $\eta = 0.05$ & $\eta = 0.1$ & $\eta = 0.25$\\
 \hline
  ${\mathcal P}_{B1}$ & 2.4175 & 3.3268 &  3.7613 &  4.9623 
 \\ 
  \hline
   ${\mathcal P}_{B2}$ & 0.0030 & 0.0208 & 0.0841 & 0.4279
 \\  
 \hline
\end{tabular}
\end{center}
\caption{Test case 2. Perturbation bounds for the noise propagation in the response matrix.}
\label{tab_test2_sin}
\end{table}

\section{Conclusions and future developments}\label{sec13}
In this work we present a first-of-its-kind approach to make the directional HiPOD approach robust when the offline data is noisy. To this aim, we modify the 
interpolation-HiPOD proposed in~\cite{lupopasini2022} by replacing interpolation techniques with ML regression, which gives rise to what we name the regression-HiPOD method.

Numerical results, although preliminary, showcase the performance of the new HiPOD approach in ideal situations, where both the noisy and clean data is available. This allows for a practical validation of the quality of the reconstruction yielded by interpolation-HiPOD against regression-HiPOD. \\
In particular, when the level of noise in the data is large, regression-HiPOD outperforms interpolation-HiPOD in terms of accuracy, by gaining up to an order in the $L^2(\Omega)$- and $H^1(\Omega)$-norm of the relative modeling error. Moreover, 
the regression-HiPOD succeeds in accurately reproducing the HiMod solution also in the presence of localized abrupt dynamics (where interpolation-HiPOD plainly fails), as in the second test case where the forcing term is localized in  narrow regions close to the boundary. 

Additionally, we provide  a new upper bound that estimates the effect of the noise level on the deformation of the subspace spanned by the left singular vectors of the response matrix. The new quantity we propose is more practical to compute and provides a more meaningful estimate with respect to some upper bounds available in the literature, as confirmed by the values in Tables~\ref{tab_test1_sin} and~\ref{tab_test2_sin}.

Concerning possible future developments, we plain to include uncertainty quantification in the analysis, to extend the directional regression-HiPOD approach to handle multiple parameters simultaneously and to model vector PDEs such as the incompressible Stokes and Navier-Stokes equations, with a view to haemodynamics modeling~\cite{PerottoBarbosa20}.



\section*{Acknowledgements}
Massimiliano Lupo Pasini thanks Dr. Vladimir Protopopescu for his valuable feedback in the preparation of this manuscript.
This work was partially supported by the program development funding resources of Oak Ridge National Laboratory, managed by UT-Battelle, LLC, for the US Department of Energy under contract DE-AC05-00OR22725.\\
Simona Perotto acknowledges the European Union’s Horizon 2020 research and innovation programme under the Marie Skłodowska-Curie Actions, grant agreement 872442
(ARIA, Accurate Roms for Industrial Applications), and the PRIN research grant n.20204LN5N5
(Advanced Polyhedral Discretisations of Heterogeneous PDEs for Multiphysics Problems).

%
\section*{Declarations}
The authors declare that they have no conflict of interest.
\\[3mm]
Data sharing not applicable to this article as no datasets were generated or analysed during the current study.

\bibliographystyle{spmpsci}      
\bibliography{sn-bibliography}

\end{document}